\documentclass{agtart_a}
\pdfoutput=1

\usepackage[all]{xypic}
\usepackage{euscript}

\title{On realizing diagrams of $\Pi$--algebras}

\author[D Blanc]{David Blanc}
\givenname{David}
\surname{Blanc}
\address{Department of Mathematics\\
University of Haifa\\\newline
31905 Haifa\\
Israel}
\email{blanc@math.haifa.ac.il}
\urladdr{}

\author[M\,W Johnson]{Mark W Johnson}
\givenname{Mark W}
\surname{Johnson}
\address{Department of Mathematics\\
Penn State Altoona\\\newline
Altoona PA 16601\\
USA}
\email{mwj3@psu.edu}
\urladdr{}

\author[J\,M Turner]{James M Turner}
\givenname{James M}
\surname{Turner}
\address{Department of Mathematics\\
Calvin College\\\newline
Grand Rapids MI 49546\\
USA}
\email{jturner@calvin.edu}
\urladdr{}

\volumenumber{6}
\issuenumber{}
\publicationyear{2006}
\papernumber{29}
\startpage{763}
\endpage{807}

\doi{}
\MR{}
\Zbl{}


\received{20 October 2005}
\revised{5 April 2006}
\accepted{5 April 2006}
\published{21 June 2006}
\publishedonline{21 June 2006}
\proposed{}
\seconded{}
\corresponding{}
\editor{}
\version{}

\arxivreference{math.AT/0604161}

\let\xysavmatrix\xymatrix
\def\xymatrix{\disablesubscriptcorrection\xysavmatrix}
\AtBeginDocument{\let\bar\wbar\let\tilde\wtilde\let\hat\what}

\swapnumbers

\makeatletter
\def\cnewtheorem#1[#2]#3{\newtheorem{#1}{#3}[section]
\expandafter\let\csname c@#1\endcsname\c@subsection}

\theoremstyle{plain}
	\cnewtheorem{thm}[subsection]{Theorem}
	\cnewtheorem{prop}[subsection]{Proposition}
	\cnewtheorem{lemma}[subsection]{Lemma}
	\cnewtheorem{cor}[subsection]{Corollary}
	\cnewtheorem{fact}[subsection]{Fact}
	\cnewtheorem{subsec}[subsection]{}
\theoremstyle{definition}
	\cnewtheorem{defn}[subsection]{Definition}
	\cnewtheorem{example}[subsection]{Example}
	\cnewtheorem{examples}[subsection]{Examples}
	\cnewtheorem{claim}[subsection]{Claim}
	\cnewtheorem{notn}[subsection]{Notation}
	\cnewtheorem{construct}[subsection]{Construction}
\theoremstyle{remark}
	\cnewtheorem{remark}[subsection]{Remark}
	\cnewtheorem{remarks}[subsection]{Remarks}
	\cnewtheorem{ack}[subsection]{Acknowledgements}

\makeautorefname{defn}{Definition}

\makeatother  

\newcommand{\mydiagram}[2][]
{\begin{equation}{#1}\vcenter{\xymatrix{#2}}\end{equation}}
\newcommand{\w}[2][]{\ensuremath{#2}{#1}}
\newcommand{\ssusp}[2]{{#1}\otimes\bS{#2}}
\newcommand{\csusp}{\Sigma_{\calC}}
\newcommand{\xra}[1]{\xrightarrow{#1}}
\newcommand{\xla}[1]{\xleftarrow{#1}}
\newcommand{\hra}{\hookrightarrow}

\newcommand{\hsp}{}
\newcommand{\hsm}{}
\newcommand{\vsm}{}

\newcommand{\lra}[1]{\langle{#1}\rangle}
\newcommand{\llra}[2]{\langle\!\langle{#1},{#2}\rangle\!\rangle}

\newcommand{\DEF}{:=}

\newcommand{\epic}{\twoheadrightarrow}
\newcommand{\equaliz}{\rightrightarrows}
\newcommand{\equalizer}{\rightrightarrows}
\newcommand{\hotimes}{\hat{\otimes}}

\newcommand{\Cok}{\operatorname{Coker}}
\newcommand{\colim}{\operatorname{colim}}

\newcommand{\gr}{\operatorname{gr}}
\newcommand{\ho}{\operatorname{ho}}
\newcommand{\holim}{\operatorname{holim}}
\newcommand{\Hom}{\operatorname{Hom}}
\newcommand{\Id}{\operatorname{Id}}
\newcommand{\incl}{\operatorname{incl}}
\newcommand{\Image}{\operatorname{Im}}
\newcommand{\Ker}{\operatorname{Ker}}
\newcommand{\map}{\operatorname{map}}
\newcommand{\Obj}{\operatorname{Obj}}
\newcommand{\op}{^{\operatorname{op}}}

\newcommand{\sk}[1]{\operatorname{sk}_{#1}}
\newcommand{\tru}[1]{\operatorname{tr}_{#1}}
\newcommand{\Alg}[1]{{#1}\text{--}{\EuScript Alg}}
\newcommand{\Gp}{{\EuScript Gp}}
\newcommand{\Pa}{$\Pi$--algebra}
\newcommand{\PAlg}{\Alg{\Pi}}
\newcommand{\PAn}{\PAlg_{n}^{n+2}}
\newcommand{\PAnk}{\PAlg_{n}^{k}}
\newcommand{\PiA}{\Pi_{\calA}}
\newcommand{\PiB}{\Pi_{\calB}}
\newcommand{\piA}{\pi_{\calA}}
\newcommand{\PAAlg}{\Alg{\PiA}}
\newcommand{\Ab}{{\operatorname{Ab}}}
\newcommand{\ab}{{\operatorname{Ab}}}
\newcommand{\APa}{$\PiA$--algebra}
\newcommand{\Ta}{$\Theta$--algebra}
\newcommand{\RM}[1]{{#1}\text{--}{\EuScript Mod}}

\newcommand{\Sets}{{\EuScript Set}_{\ast}}
\newcommand{\Spec}{{\EuScript Spec}}
\newcommand{\Top}{\calT}
\newcommand{\DD}{\mathbb D}
\newcommand{\LL}{\mathbb L}
\newcommand{\NN}{\mathbb N}
\newcommand{\QQ}{\mathbb Q}
\newcommand{\RR}{\mathbb R}
\newcommand{\mS}{\mathbb S}
\newcommand{\ZZ}{\mathbb Z}
\newcommand{\fG}{\mathfrak{G}}
\newcommand{\calA}{{\mathcal A}}
\newcommand{\calB}{{\mathcal B}}
\newcommand{\calC}{{\mathcal C}}
\newcommand{\sC}{s\calC}
\newcommand{\CA}{\sC_{\calA}}
\newcommand{\hC}{\ho\calC}

\newcommand{\calG}{{\mathcal G}}

\newcommand{\MC}{\calC(\to)}
\newcommand{\calO}{{\mathcal O}}
\newcommand{\calS}{{\mathcal S}}
\newcommand{\calSs}{\calS_{\ast}}
\newcommand{\calT}{{\mathcal T}}
\newcommand{\calTs}{\calT_{\ast}}

\numberwithin{equation}{section}
\newcommand{\Delt}[1]{\Delta[{#1}]}
\newcommand{\dD}[1]{\partial\Delt{#1}}
\newcommand{\bV}{\mathbf{V}}
\newcommand{\bW}{\mathbf{W}}
\newcommand{\bX}{\mathbf{X}}
\newcommand{\bY}{\mathbf{Y}}
\newcommand{\bZ}{\mathbf{Z}}
\newcommand{\q}[1]{^{({#1})}}
\newcommand{\cons}[1]{c({#1})_{\bullet}}

\newcommand{\bp}[1]{\mathbf{p}\q{#1}}

\newcommand{\ruq}[1]{\rho\q{#1}}

\newcommand{\ts}[1]{\tilde{s}_{#1}}

\newcommand{\pis}{\pi_{\ast}}
\newcommand{\pinat}{\pi^{\natural}}
\newcommand{\pin}[2]{\pinat_{{#1}}{#2}}
\newcommand{\HL}[3]{H^{#1}_{\Lambda}({#2};\,{#3})}
\newcommand{\exal}[3]{\operatorname{exal}_{#1}({#2};\,{#3})}
\newcommand{\Bd}{B_{\bullet}}
\newcommand{\Ed}{E_{\bullet}}
\newcommand{\Fd}{F_{\bullet}}
\newcommand{\Gd}{G_{\bullet}}

\newcommand{\Kd}{{K_{\bullet}}}
\newcommand{\Ld}{L_{\bullet}}
\newcommand{\Md}{M_{\bullet}}

\newcommand{\Nd}{N_{\bullet}}
\newcommand{\Pd}{P_{\bullet}}
\newcommand{\Vd}{V_{\bullet}}
\newcommand{\Wd}{W_{\bullet}}
\newcommand{\Xd}{X_{\bullet}}
\newcommand{\Xn}[1]{X\lra{{#1}}_{\bullet}}
\newcommand{\fn}[1]{f\lra{{#1}}}
\newcommand{\Yd}{Y_{\bullet}}
\newcommand{\Yn}[1]{Y\lra{{#1}}_{\bullet}}
\newcommand{\gS}[1]{{\EuScript S}^{#1}}
\newcommand{\bS}[1]{\mathbf{S}^{#1}}
\newcommand{\be}[1]{\mathbf{e}^{#1}}
\newcommand{\eV}{{\EuScript V}}
\newcommand{\eVd}{\eV_{\bullet}}
\newcommand{\eW}{{\EuScript W}}
\newcommand{\eWd}{\eW_{\bullet}}
\newcommand{\bVd}{\bV_{\bullet}}

\newcommand{\Ts}{T_{\ast}}
\newcommand{\As}{\Lambda}
\newcommand{\OA}[1]{\Omega^{#1}\Lambda}
\newcommand{\BL}{B\Lambda}
\newcommand{\BCL}[1]{B_{{#1}}\Lambda}
\newcommand{\EM}[3]{E\sp{{#1}}({#2},{#3})}
\newcommand{\EC}[3]{E\sb{{#1}}({#2},{#3})}
\newcommand{\EL}[2]{\EM{\Lambda}{#1}{#2}}
\newcommand{\ECL}[3]{E\sp{\Lambda}\sb{{#1}}({#2},{#3})}
\newcommand{\tB}{\tilde{B}}
\newcommand{\tBL}{\tB\Lambda}
\newcommand{\tE}[3]{\tilde{E}\sp{{#1}}({#2},{#3})}
\newcommand{\tEL}[2]{\tE{\Lambda}{#1}{#2}}
\newcommand{\tP}[1]{\tilde{P}_{#1}}
\newcommand{\tk}[1]{\tilde{k}_{#1}}
\newcommand{\bB}{\mathbf{B}}
\newcommand{\bBL}{\bB\Lambda}
\newcommand{\bE}[3]{\mathbf{E}\sp{#1}({#2},{#3})}
\newcommand{\bEL}[2]{\bE{\Lambda}{#1}{#2}}
\newcommand{\bP}[1]{\mathbf{P}_{#1}}
\newcommand{\bk}[1]{\mathbf{k}_{#1}}
\begin{document}

\begin{htmlabstract}
Given a diagram of &Pi;&ndash;algebras (graded groups equipped with an action
of the primary homotopy operations), we ask whether it can be realized
as the homotopy groups of a diagram of spaces. The answer given here
is in the form of an obstruction theory, of somewhat wider application,
formulated in terms of <em>generalized &Pi;&ndash;algebras</em>. This extends a
program begun in [J Pure Appl Alg. 103 (1995) 167-188] and [Topology
43 (2004) 857-892] to study the realization of a single &Pi;&ndash;algebra.
In particular, we explicitly analyze the simple case of a single map,
and provide a detailed example, illustrating the connections to higher
homotopy operations.
\end{htmlabstract}

\begin{webabstract}
Given a diagram of $\Pi$--algebras (graded groups equipped with an action
of the primary homotopy operations), we ask whether it can be realized
as the homotopy groups of a diagram of spaces. The answer given here
is in the form of an obstruction theory, of somewhat wider application,
formulated in terms of \emph{generalized $\Pi$--algebras}. This extends a
program begun in [J.\ Pure Appl.\ Alg. 103 (1995) 167-188] and [Topology
43 (2004) 857-892] to study the realization of a single $\Pi$--algebra.
In particular, we explicitly analyze the simple case of a single map,
and provide a detailed example, illustrating the connections to higher
homotopy operations.
\end{webabstract}

\begin{asciiabstract}
Given a diagram of Pi-algebras (graded groups equipped with an action
of the primary homotopy operations), we ask whether it can be realized
as the homotopy groups of a diagram of spaces. The answer given here
is in the form of an obstruction theory, of somewhat wider application,
formulated in terms of generalized Pi-algebras. This extends a program
begun in [J.\ Pure Appl.\ Alg. 103 (1995) 167-188] and [Topology 43
(2004) 857-892] to study the realization of a single Pi-algebra.
In particular, we explicitly analyze the simple case of a single map,
and provide a detailed example, illustrating the connections to higher
homotopy operations.
\end{asciiabstract}

\begin{abstract}
Given a diagram of \Pa s (graded groups equipped with an action of
the primary homotopy operations), we ask whether it can be realized as
the homotopy groups of a diagram of spaces. The answer given here is
in the form of an obstruction theory, of somewhat wider application,
formulated in terms of \emph{generalized \Pa s}. This extends a program
begun by Dwyer, Kan, Stover, Blanc and Goerss \cite{DKStoE,BDGoerR} to
study the realization of a single \Pa.  In particular, we explicitly
analyze the simple case of a single map, and provide a detailed
example, illustrating the connections to higher homotopy operations.
\end{abstract}

\maketitle

%
%
\section{Introduction}
\label{cint}

A recurring problem in algebraic topology is the rectification of
homotopy-commutative diagrams:	given a diagram
\w{F\co\DD\to\ho\calTs} (that is, a functor from a small category to the
homotopy category of topological spaces), we ask whether $F$ lifts to
\w[,]{\hat{F}\co\DD\to\calTs} and if so, in how many ways.

Such questions arise naturally in determining if a given
$H$--space is a loop space; in defining Steenrod operations; in
analyzing structured ring spectra; and so on. Our goal here is to
present an obstruction-theoretic approach to an algebraic version of
this question.

\subsection{Diagrams of \Pa s}\label{spad}
Recall that a \Pa\ is a graded group equipped with an action of the
primary homotopy operations (Whitehead products and compositions),
modeled on the homotopy groups of a space (see \fullref{ccpa} below).
In \cite{DKStoE,DKStoB}, Dwyer, Kan, and Stover set out to construct
an obstruction theory for realizing a given \Pa\ $\Lambda$ as
\w[,]{\Lambda\cong\pis X} for some space $X$. The program was completed
by Blanc, Dwyer and Goerss in \cite{BDGoerR}, using methods developed
by Dwyer and Kan in a series of papers which established a general
obstruction theory for rectifying homotopy-commutative diagrams (see
the work of Dwyer, Kan and Smith \cite{DKanC,DKanO,DKanR,DKanCR,DKSmiR}).
Our goal here is to extend this program to address the following: %

\subsection{Diagram realization question}
Can a given diagram of \Pa s \w{\Lambda\co\DD\to\PAlg} be \emph{realized} --
that is, lifted to a diagram of spaces \w{\hat{\Lambda}\co\DD\to\calTs} with
\w[?]{\pis\circ\hat{\Lambda}=\Lambda\vsm}

The answer we provide is in the form an obstruction theory: we
inductively define a sequence of cohomology classes
\w[,]{k_{n}\in H^{n+2}(\Lambda;\Omega^{n}\Lambda)} and
show that $\Lambda$ is realizable precisely when \w{k_{n}=0}
for all $n$. The case of a single \Pa\	was treated in \cite{BDGoerR},
and the extension to our context is straightforward. However, the
description there was in terms of moduli spaces, and it seems
worthwhile making obstruction theory explicit.
A further generalization of this theory appears in \cite{BlaCH}, but
it is not easy to extract from it the simpler version needed here.

\subsection{Generalized \Pa s}\label{sgpa}
In fact, it turns out that this approach may be carried out somewhat more
generally, for any $E^{2}$--model category \w{\sC} (see \fullref{cemc}),
once we have chosen a set $\calA$ of homotopy cogroup objects
in $\calC$ to play the role of the spheres \w{\{\bS{n}\}_{n=1}^{\infty}}
in \w[.]{\calTs}

Note that a \Pa\ can be thought of as a product-preserving functor
\w[,]{T\co \Pi\op\to\Sets} where $\Pi$ is the subcategory of finite wedges of
spheres in \w[.]{\ho\calTs}  Similarly defining \w{\PiA\subseteq\ho\calC}
for any $\calA$ as above, we define a \emph{\APa} to be a
product-preserving functor \w[.]{\PiA\op\to\Sets}

For example, a map \w{\phi\co \Gamma\to\Lambda} of ordinary
\Pa s corresponds to a diagram in \w[,]{(\PAAlg)^{\DD}}
where $\DD$ has two objects and a single non-identity map \w[.]{0\to 1}
Setting
\[
\calA:=\{\bS{n}\stackrel{\Id}{\to}\bS{n},\hsm \ast\to\bS{n}\}_{n\in\NN},
\]
we can think of $\phi$ as a generalized \APa. The realization question
for diagrams of \Pa s is thus a special case of the the following\vsm:

\subsection{General Realization Question}
Given a model category $\calC$ with set of models $\calA$, when is a
\APa\ $\Lambda$ realizable in $\calC$? That is, is there an
\w{X\in\calC} such that \w{\piA X\cong \Lambda} (where \w{\piA X} is
defined by
\w[?)]{A \mapsto [A,X]_{\calC}\vsm}

Again, this is not meant to be a gratuitous exercise in generalization:
it allows us to answer in a systematic way the same question for
(diagrams of) localized or $n$--connected spaces, spectra, $n$--types,
and so on.

\subsection{Notation and conventions}\label{snac}
$\calT$ will denote the category of topological spaces, and \w{\calTs}
that of pointed connected topological spaces. By a
\emph{space} we shall always mean an object in \w[.]{\calTs}

The category of groups is denoted by $\Gp$, and that of pointed
sets by \w[.]{\Sets} For any category $\calC$, \w{\gr_{\calA}\calC}
denotes the
category of  $\calA$--graded objects over $\calC$  (that is, the category
\w{\calC^{\calA}} of diagrams indexed by the discrete category $\calA$),
and \w{\sC} that of simplicial objects over $\calC$. The category of
simplical
sets will be denoted by $\calS$, that of pointed connected
simplicial sets by \w[,]{\calSs} and that of simplicial groups by $\calG$.
For any \w[,]{Z\in\calC} write \w{\cons{Z}} for the constant
simplicial object determined by $Z$.

The suspension in a model category $\calC$ will denote the usual
pushout of the inclusions in two cones (that is, factorizations of the
final map as a cofibration followed by an acyclic fibration),
following Quillen \cite[Section~I.2]{QuiH}.  This operation will be
indicated by \w{\csusp} henceforth.

\begin{defn}\label{dskel}
The category of simplicial objects \w{X_{0},\dots,X_{n}} truncated at the
$n$th dimension will be denoted by \w[.]{s_{n}{\calC}} If $\calC$
has enough
colimits, the obvious truncation functor \w{\tru{n}\co \sC\to s_{n}\calC} has
a left adjoint $\rho_{n}\co s_{n}\calC\to\sC$, and the composite
\w{\sk{n}\DEF \rho_{n}\circ\tru{n}\co \sC\to\sC} is called the
$n$--\emph{skeleton} functor. Thus \w{\sk{n}\Xd} is ``freely
generated'' as a simplicial object by \w[.]{X_{0},\ldots, X_{n}}
\end{defn}

\begin{defn}\label{dsimpst}
Let \w{\Delt{n}} denote the standard $n$--simplex in $\calS$, generated by
\w[,]{\sigma_{n}\in \Delt{n}_{n}} with boundary \w{\dD{n}} (the
sub-object generated by \w{d_{i}\sigma_{n}} for \w[).]{0\leq i\leq n}
Similarly, the $k$th-horn \w{\Lambda^{k}[n]} is the sub-object generated
by
\w{d_{i}\sigma_{n}} for \w[.]{i\neq k} The simplicial $n$--sphere is
\w[.]{\bS{n}\DEF\Delt{n}/\dD{n}}

If $\calC$ has enough colimits, for \w{M\in\calSs} and \w[,]{X\in\calC}
we define \w{X\hotimes M\in\sC} by
\w[,]{(X\hotimes M)_{n}:=\coprod_{m\in M_{n}} X}
with face and degeneracy maps induced from those of $M$.
For \w[,]{Y\in\sC} define \w{Y\otimes M\in\sC} by
\w[.]{(Y\otimes M)_{n}:= \coprod_{m\in M_{n}} Y_{m}}
The simplicial suspension functor \w{\ssusp{-}{n}} (on \w{\sC}) is
defined by
\w[.]{\ssusp{Y}{n}:=Y\otimes(\Delt{n}/\partial\Delt{n})}
\end{defn}

The main result of this paper is an obstruction theory for dealing
with the general realization question, expressed in the following:

%
%
\begin{thm}[Theorems~\ref{tfour} and~\ref{tclass}]\label{tzero}
A \APa\ $\Lambda$ can be realized in $\calC$ if and only if an
inductively-defined sequence of cohomology classes in
\w{\HL{n+3}{\Lambda}{\OA{n+1}}} all vanish. The different realizations
(if any) are classified (up to homotopy) by elements of
\w[.]{\HL{n+2}{\Lambda}{\OA{n+1}}}
\end{thm}

\subsection{Higher homotopy operations}\label{shho}

Higher order homotopy operations appear as obstructions to rectifying
homotopy commutative diagrams, so, as one might expect, they tie in
with our approach (in more than one way).
One of the original motivations for this paper was to try to understand
the intriguing relationship between the diagram realization question,
framed in
the algebraic language of \Pa s and cohomology, and the motivating
topological problem of rectifying homotopy commutative diagrams.
A general answer is still beyond us (but see \fullref{rtcoh} below).
We shall, however, show how this connection appears in a specific example,
which we will be using as a leitmotif to illustrate various
constructions throughout this paper.

\begin{defn}\label{dtoda}
Given a homotopy commutative diagram
%
\begin{equation}
\label{eqtwentyeight}
\xymatrix{\bW \ar[r]^{f} \ar@/^{1.5pc}/[rr]^{\ast} & \bX \ar[r]^{g}
\ar@/_{1pc}/[rr]_{\ast} & \bY \ar[r]^{h} & \bZ}
\end{equation}
\noindent the \emph{Toda bracket} \w{\lra{f,g,h}\subseteq[\Sigma\bW,\bZ]}
is the set of all homotopy classes which are pushout maps $k$ in the
diagram
\begin{equation}{\label{epush}}
\vcenter{\xymatrix@C=40pt{
\ar @{} [dr] |>>>>>{\framebox{\scriptsize{PO}}}
\bW  \ar[r]^-{i_{1}} \ar[d]^-{i_{2}} & C\bW \ar[d]
\ar@/^/[ddr]^{G\circ Cf}  & \\
C\bW \ar[r] \ar@/_/[drr]_{h\circ F} & \Sigma\bW \ar@{.>}[dr]^{k}& \\
& & \bZ
}}
\end{equation}
\noindent where \w{G\co h\circ g\sim \ast} and \w{F\co g\circ f\sim \ast}
are any
nullhomotopies.
\end{defn}

Note that \w{\lra{f,g,h}} is the obstruction to rectifying the
homotopy commutative diagram \eqref{eqtwentyeight}, in the sense that
it vanishes (that is, contains the null class) if and only if
\eqref{eqtwentyeight} can be rectified (that is, realized by a strictly
commutative diagram, with the null maps represented by actual zero maps).

\begin{example}\label{etoda}
Recall that in the stable range
%
\begin{equation}\label{eqtwentyseven}
\pi_{i}\bS{k}\cong
\begin{cases}
\ZZ\lra{\iota} & \text{for } i=k\\
(\ZZ/2)\lra{\eta} & \text{for } i=k+1\\
(\ZZ/4)\lra{\eta^{2}} & \text{for } i=k+2\\
(\ZZ/24)\lra{\nu} & \text{for } i=k+3\\
0 & \text{for } i=k+4, k+5\\
\end{cases}
\end{equation}
\noindent where \w{\eta^{3}=12\nu} (cf Toda \cite[14.1]{TodaC}).
Thus, for $k\geq 3$, the sequence
$$
\xymatrix@R=25pt{
\bS{k+2}\ar[r]^{\eta} & \bS{k+1}\ar[r]^{2} & \bS{k+1}\ar[r]^{\eta} &\bS{k}
}
$$
\noindent is an instance of \eqref{eqtwentyeight}, with the
corresponding Toda bracket
%
\begin{equation}\label{eqtwentynine}
\lra{\eta,2,\eta}=\{\nu,\nu+\eta^{3}\}=\{\pm\nu\}\subseteq\pi_{k+3}\bS{k}.
\end{equation}
\noindent (See Toda \cite[(5.4)]{TodaC}).
\end{example}

\begin{remark}\label{rtcoh}
Given a homotopy-commutative diagram \w{F\co \DD\to\ho\calTs} of
topological spaces (for most reasonable indexing diagrams $\DD$),
a suitable higher homotopy operation appears as the obstruction to
rectifying $F$ (that is, lifting it to \w[).]{\calTs} However, in many
applications all spaces in the diagram (except perhaps
\w[,]{F(\ast)} where $\ast$ is terminal in $\DD$) are
(wedges of) spheres -- as in \fullref{etoda}.

In this case we can replace $F$ by the corresponding diagram of \Pa s
\w{\pis\circ F\co \DD\to\PAlg} with no loss of generality (beyond the
choice of realization for \w[),]{\pis F(\ast)} and any obstruction to
realizing \w{\pis\circ F} is in particular an obstruction to
rectifying $F$. Thus \fullref{tzero} provides a way to describe
many higher homotopy operations algebraically, in terms of suitable
cohomology classes. We hope to pursue this point further in a future
paper.
\end{remark}

\subsection{Organization}

In \fullref{ccpa} we define our objects of study,
\APa s and some related algebraic concepts.
\fullref{crmc} begins a detailed analysis of
resolution model category structures on \w[,]{\sC} and their basic
properties, giving several important examples. \fullref{cemc}
defines $E^{2}$--\emph{model categories}, which are a special kind of
resolution
model category provided with additional structures, such as
Eilenberg--Mac\,Lane objects and Postnikov towers.
The motivating examples of diagram categories of spaces, as well as
the main algebraic categories, are all instances of this.  In fact, we
show that any diagram category on an $E^{2}$--model category is another,
which provides a broad class of examples.

In \fullref{ccoh}, we define the cohomology theory associated
to an $E^{2}$--model category structure and describe some of its
basic properties.  We illustrate this for the simplest example of a
diagram category, namely an arrow category, and show how the
cohomology of an arrow relates to that of the source and target objects.

The technical heart of the paper is the obstruction
theory for dealing with the general realization question, which appears
in \fullref{crpia}.  As expected, we induct up the construction of
the Postnikov tower of our (putative) simplicial object expected to
yield a realization of $\Lambda$.
\fullref{crfpi} provides a more explicit description of the
single map case, illustrating it with a detailed example.

\begin{ack}
We would like to thank the referee for his or her comments.
The third author was supported by NSF grant DMS-0206647 and a
   Calvin Research Fellowship (SDG).
\end{ack}

%
%
\section{\APa s}
\label{ccpa}

The functor \w{X\mapsto\pis X} is corepresented by spheres in the
homotopy category of spaces. If we want to include the group
structures, Whitehead products,
and $\pi_1$--actions as well, we expand the domain category
(choices of the argument $?$ for \w[)]{[?,X]} to finite wedges of spheres,
and require that wedges be sent to products.
This definition extends to other model categories, using the relevant
properties of spheres:

\begin{defn}\label{models}
Let $\calC$ be a cofibrantly generated pointed model category which is
  \emph{right proper} -- that is, the pullback of a weak
  equivalence along a fibration is a weak equivalence. A collection of
  \emph{models} for $\calC$ is a set $\calA$ of cofibrant
  homotopy cogroup objects in $\calC$, closed under suspension in
  $\calC$ (denoted by \w[).]{\csusp}
\end{defn}

\begin{defn}\label{dcpa}
Given a model category $\calC$ as above and a set $\calA$ of models
for $\calC$, let \w{\PiA} denote the full subcategory of \w{\hC}
consisting of fibrant and cofibrant objects weakly equivalent to
finite coproducts of objects from $\calA$ (which become
products in \w[).]{\PiA\op} A \emph{\APa} is defined to be a
product-preserving functor \w[,]{\PiA\op\to\Sets}
and the category of \APa s (and natural transformations) will be
denoted by \w[.]{\PAAlg}
\end{defn}

Since the suspension operator in $\calC$ preserves the class of cofibrant
homotopy cogroup objects, in many of our examples $\calA$ is generated
under \w{\csusp} by a much smaller set.  For example, the set of
spheres used to define ordinary \Pa s is generated by the circle
\w[.]{\bS{1}}

\begin{example}\label{rcpa}
The canonical example of a \APa\ is a \emph{realizable} \APa --
that is,
one given by \w{[?,X]_{\calC}} for some \w[.]{X \in \calC}  This will
be referred to as the \emph{homotopy \APa} of $X$; it defines a functor
\w[.]{\piA\co \hC\to\PAAlg}
\end{example}

\begin{remark}\label{cpa}
When \w{\calC=\Gp} is the category of groups, and
\w[,]{\calA=\{\ZZ\}} the category of \APa s is equivalent to \w{\Gp}
itself.  In \fullref{ermc}(f), we allow for a range of
universal algebras as examples for $\calC$. As noted by Quillen
\cite[Section~II]{QuiH}, there is an (unique) object \w{D \in \calC} such
that, for
\w[,]{\calA =\{D\}} the category \w{\PAAlg} is equivalent to $\calC$.

On the other hand, in the resulting resolution model category
\w{\calG=\sC} with \w[,]{\calA=\{\ZZ\}} (under the constant embedding of
$\calC$ in \w[),]{\sC} the category \w[,]{\PiA} consisting of all
suspensions
of $\ZZ$ and coproducts thereof, is just the $\calG$--version of the
collection
of all wedges of spheres (in \w[),]{\Ts} so \w{\PAAlg} is the original
category of \Pa s (cf Stover \cite[Section~2]{StoV}). See \fullref{etoda} and
\fullref{sspace} for examples of such \Pa s.
\end{remark}

\subsection{The free functor}\label{scpia}

There is a forgetful functor \w{\calO\co \PAAlg\to\gr_{\calA}\Sets}
to the category of $\calA$--graded pointed sets, with left adjoint
\w[.]{F\co \gr_{\calA}\Sets\to\PAAlg} We call \w{F(W)}
the \emph{free} \APa\ \emph{generated by} \w[.]{W\in\gr_{\calA}\Sets}
Thus \w{\PAAlg} is an FP-sketchable variety of universal algebras as
in \fullref{ermc}(f), sketched by the $\fG$--theory
\w[.]{\Theta:=\PiA} In particular, \w{\PAAlg}  is complete and
cocomplete (see Ad\'amek and Rosick\'y \cite[Section~1]{ARosiL}).

\subsection{Products and coproducts}\label{ecoprod}

We now describe a variety of constructions which will be used at
various points later.  Given two \APa s $\Lambda$ and $\Gamma$	over a
fixed \APa\ $B$, we define their \emph{fibered product}
\w{\Lambda\times_{B}\Gamma} in \w{\PAAlg/B} by declaring its value on
an object \w{U\in \PiA} to be the set-theoretic pullback
\begin{equation}\label{epull}
\vcenter{\xymatrix@C=10pt{
      {(\Lambda\times_{B} \Gamma)(U)} \ar[d] \ar[rr] &&
	 {\prod_{\beta}(\Lambda(U)\times_{\prod_{\gamma}B(U_{\gamma})}
	 \Gamma (U_{\beta}))} \ar[d]^{\hbox{\footnotesize$\sim$}} \\
      {\prod_{\alpha}(\Lambda(U_{\alpha})\times_{\prod_{\gamma}B(U_{\gamma})}
      \Gamma (U))} \ar[rr]^-{\hbox{\footnotesize$\sim$}} &&
	{\prod_{\alpha}\prod_{\beta}(\Lambda(U_{\alpha})
	\times_{\prod_{\gamma}B(U_{\gamma})} \Gamma (U_{\beta}))}
}}
\end{equation}
whenever \w{U = \coprod_{\alpha}U_{\alpha}} for
\w[.]{U_{\alpha}\in\PiA}

Similarily, the \emph{coproduct} \w{\Lambda_{0} \amalg \Lambda_{1}} of two
\APa s \w{\Lambda_{0}} and \w{\Lambda_{1}} may be characterized
explicitly by first setting
\w[,]{\Lambda_{0}\amalg\Lambda_{1}:= F(W_{0}\vee W_{1})} if
\w{\Lambda_{0}=F(W_{0})} and \w{\Lambda_{1} = F(W_{1})} are free;
and, more
generally, as the natural group quotient
$$
(F\calO \Lambda_{0}\amalg F\calO \Lambda_{1})/I
$$
where $I$ is the smallest ideal containing the kernels \w{K_{i}} of
\w{F\calO \Lambda_{i} \to \Lambda_{i}} for \w[.]{i = 0,1} Note there is
also a
coequalizer in \w{\PAAlg}
$$
(F\calO)^{2}\Lambda_{0}\amalg (F\calO)^{2}\Lambda_{1}\equalizer
(F\calO)\Lambda_{0}\amalg
(F\calO)\Lambda_{1}\to\Lambda_{0}\amalg\Lambda_{1}
$$
induced by the two adjunction maps \w{F\calO\to\Id} and \w[.]{\Id\to\calO
F}

\begin{defn}\label{dideal}
An \emph{ideal} in a \APa\ $\Lambda$ is a sub-\APa\
\w[,]{I\subset\Lambda} such that for any \w[,]{U\in\PiA} the top arrow
in the commuting diagram
\mydiagram{
{\Lambda(U)\times I(U)} \ar[r] \ar[d]
       & {I(U)} \ar[d] \\
{\Lambda(U)\times \Lambda (U)} \ar[r] & {\Lambda(U)}
}
\noindent exists. (Uniqueness follows from injectivity of
\w[).]{I(U)\to\Lambda(U)} For example, the kernel
\w{\Ker(f):=\ast\times_{\Gamma}\Lambda} of a map of \APa s
\w{f\co \Lambda \to \Gamma} is an ideal.
\end{defn}

\begin{defn}\label{dlalg}
For a fixed \APa\  $\Lambda$, a \emph{$\Lambda$--\APa} is a map of \APa s
\w[.]{i\co  \Lambda \to \Gamma} In particular, given
\w[,]{W\in\gr_{\calA}\Sets}
the \emph{free} $\Lambda$--\APa\ on $W$ is
defined by \w[.]{F_{\Lambda}(W) := F(W)\amalg \Lambda} Similarly, we can
define the
\emph{$\Lambda$--coproduct} \w{\Gamma_{1}\amalg_{\Lambda}\Gamma_{2}}
of two $\Lambda$--\APa s \w{\Gamma_{1}} and \w{\Gamma_{2}} as a
coequalizer in \w{\PAAlg}
$$
 \Lambda \equaliz
\Gamma_{1}\amalg\Gamma_{2}\to\Gamma_{1}\amalg_{\Lambda}\Gamma_{2}
$$
\noindent where the left pair of maps is defined using the maps to
left/right factors  \w{\Lambda \equaliz \Lambda\amalg\Lambda} together
with the coproduct of the $\Lambda$--algebra structure maps for
\w[,]{\Gamma_{i}} \w[.]{i=1,2}

Given an ideal \w[,]{I\subseteq\Lambda} the \emph{quotient} \APa\ of
$\Lambda$ by $I$ is then defined: \w[.]{\Lambda/I :=
\ast\amalg_{I}\Lambda}
\end{defn}

\begin{defn}\label{eloop}
If $\Lambda$ is a \APa, we define the \emph{loop} \APa\
\w{\Omega\Lambda} by \w[,]{\Omega\Lambda(U):=\Lambda(\csusp U)}
where \w{\csusp U} is the suspension of $U$ in $\calC$.
\end{defn}

\subsection{Abelian \APa s}\label{sapia}

An abelian group object $M$ in \w{\PAAlg} is called an
\emph{abelian \APa} -- that is, if \w{\Hom_{\PAAlg}(B,M)} has a
natural abelian group structure for any $B$.
Note that the structure is induced by the underlying $\calA$--graded
group structure in \w[,]{\PAAlg} so in particular \w{\calO M} is an
$\calA$--graded \emph{abelian} group.

Denote by \w{\ab(\PAAlg)} the subcategory of abelian \APa s. The
inclusion functor \w{\ab(\PAAlg)\to\PAAlg} has a left adjoint
\w[,]{\Ab} called the \emph{abelianization functor}, defined for
\w{\Lambda = F(W)} by
$$
(\Ab{(F(W))})(A) := \oplus_{W_{A}} \Ab{(\piA(A))}.
$$
\noindent For general $\Lambda$, define \w{\Ab(\Lambda)} to be the
coequalizer in \w{\PAAlg}
$$
\Ab((F\calO)^{2}\Lambda)\equaliz
\Ab((F\calO)\Lambda)\to\Ab(\Lambda).
$$
Note that the composite \w{\Ab\circ F\co \gr_{\calA}\Sets\to\ab(\PAAlg)}
is left adjoint to the forgetful functor, so it is the \emph{free
abelian \APa} functor. From the adjointness we get a natural
abelianization map \w{\rho\co \Lambda\to\Ab(\Lambda)}
and we define the ideal \w{W(\Lambda)\subseteq\Lambda} as
\w[.]{\Ker(\rho)}

Then \w{W(\Lambda)} may be viewed as the \emph{ideal of primary operations
acting on elements of $\Lambda$}, and we have
\w[.]{\Lambda/W(\Lambda)\cong\Ab(\Lambda)}

\subsection{Modules}\label{smod}

For a fixed \APa\  $\Lambda$, a \emph{module over $\Lambda$} is
an abelian group object \w{p\co M\to\Lambda} in the over-category
\w[.]{\PAAlg/\Lambda}  This means that it is endowed with maps
$$
m\co M\times_{\Lambda}M \to M \text{ and } i\co M\to M
$$
\noindent in \w[,]{\PAAlg/\Lambda} as well as a section
\w{s\co \Lambda\to M} for $p$ (which represents the unit element in the
abelian group \w[).]{\Hom_{\Lambda}(\Lambda ,M)}
The category of modules over $\Lambda$ is denoted by \w[.]{\RM{\Lambda}}

Moreover, given a map of \APa s \w[,]{\Lambda \to \Gamma} the
associated restriction functor \w{\RM{\Gamma} \to \RM{\Lambda}} has a
left adjoint, which we denote by \w[.]{(-)\ast_{\Lambda}\Gamma}

Note that \w{K\DEF\Ker(p)} is itself an abelian \APa, as we can see by
mapping \w{0\co X\to\Lambda} to \w{p\co M\to\Lambda} in \w{\PAAlg/\Lambda}
for any \APa\ $X$, so we have a split exact sequence of \APa s
\begin{equation}
\label{algext}
\xymatrix{0 \rto & K \rto & M\ar[r] & \ar@/_1pc/[l] \Lambda \rto & 0,}
\end{equation}
\noindent and in particular \w{\calO M=\calO\Lambda\ltimes\calO K} is
a semi-direct product of groups.

However, $K$ is not just an abelian \APa; it also has an action of
$\Lambda$ on it, determined by an \emph{action map}
$$
\phi_{f}\co \Lambda(U)\ltimes K(U)\to K(V)
$$
for each \w{f\co V\to U} in \w[,]{\PiA}, subject to the requirements that:
\begin{enumerate}
\item The composite \w{K(U)\to\Lambda(U)\times K(U)
      \stackrel{\phi_{f}}{\longrightarrow} K(V)} is equal to \w[;]{K(f)}

\item For \w{g\co W\to V} in \w[,]{\PiA} the action map \w{\phi_{f\circ g}}
     equals the composite
$$
      \Lambda(U)\times K(U)
      \stackrel{\Delta\times\Id}{\longrightarrow}
      \Lambda(U)\times(\Lambda(U)\times K(U))
      \stackrel{\Lambda(f)\times \phi_{f}}{\longrightarrow}
      \Lambda(V)\times K(V) \xra{\phi_{g}}K(W)
$$
\end{enumerate}

We sometimes say that $K$ itself, endowed with this action of
$\Lambda$, is a $\Lambda$--\emph{module} (which corresponds to the
traditional description of an $R$--module, for a ring $R$), and write
\w[.]{M=\Lambda\ltimes K}

Note that  \w[,]{\Ab\circ F_{\Lambda} \cong\Lambda\ltimes(\Ab\circ F)} so
\w{\Ab\circ F_{\Lambda}\co \gr_{\calA}\Sets\to\RM{\Lambda}} can be viewed
as the free $\Lambda$--module functor.

\begin{remark}\label{rmod}
When \w[,]{\PAAlg=\PAlg} a $\Lambda$--module $K$ is simply
an abelian \Pa, equipped with mappings
\w[,]{\llra{\,}{\,}\co \Lambda_{p}\times K_{q}\to K_{p+q}}
commuting with compositions, such that for each \w[,]{q>0}
\w{\alpha\circ x:=\llra{\alpha}{x}-x} defines an action of
\w{\Lambda_{0}} on \w[,]{K_{q}} satisfying
\w[,]{\lra{b,a}\circ(a\circ x)=-{\lra{a,b}}\circ x -\llra{a}{x}}
while for \w[,]{p>0} \w{\llra{\,}{\,}\co \Lambda_{p}\times K_{q}\to K_{p+q}}
is bilinear, and satisfies
\[
\llra{\alpha}{\llra{\beta}{x}}=
\llra{\lra{\alpha,\beta}}{x}+(-1)^{pq}\llra{\beta}{\llra{\alpha}{x}}.
\]
\end{remark}

\begin{example}\label{eloopp}
For a \APa\ $\Lambda$, define the \APa\ \w{\Omega_{+}\Lambda} by
$$
\Omega_{+}\Lambda(A):=\Lambda\Bigl((\csusp A) \bigvee A\Bigr)\,.
$$
There is then a split exact sequence
\begin{equation}
\xymatrix{\ast \rto & \Omega\Lambda \rto & \Omega_{+}\Lambda \rto &
\ar@/_1pc/[l] \Lambda \rto & \ast,}
\end{equation}
\noindent which gives \w{\Omega_{+}\Lambda} the structure of a module
over $\Lambda$.
\end{example}

\begin{example}
The \emph{fold} map \w{\nabla \co  \Lambda\amalg\Lambda \to \Lambda}
possesses two sections. Let \w[.]{K := \Ker(\nabla)} Define the
\emph{K\"{a}hler differentials} of $\Lambda$ by
\w[.]{\Omega_{\Lambda}:=\Ab(K)}
Then the split exact sequence
\begin{equation}
\label{Kaehler}
\xymatrix{\ast\rto & \Omega_{\Lambda} \rto &
\mbox{\w{\Omega_{\Lambda}\amalg_{K}(\Lambda\amalg\Lambda)}}
\rto & \ar@/_1pc/[l] \Lambda \rto & \ast}
\end{equation}
\noindent gives \w{\Omega_{\Lambda}} the structure of a $\Lambda$--module.
\end{example}

We will see in \fullref{rcoh} that the K\"{a}hler
differentials are closely related to our cohomology theories\vsm.

Our key examples of modules come in \fullref{pfive},
where we will see that for \w[,]{n>0} the natural homotopy groups
 \w{\pinat_{n} \Yd} (see \fullref{rcpa}) and their loop algebras
are modules over \w[.]{\pinat_{0} \Yd}

\begin{remark}\label{rtwocat}
We have in view two types of categories for $\calC$ here:
one type are ``algebraic'' categories, such as \w{\Gp} and
\w[,]{\PAAlg} in which the model category structures are trivial (in
the sense that the only weak equivalences are isomorphisms), so the
associated realization question is also trivial.

The other type is ``topological'' -- for example,
$\calG$ or \w[.]{\Ts} Here the associated algebraic invariants, such
as homotopy groups, give rise to meaningful realization questions; and
the associated simplicial categories possess nontrivial resolution model
category structures, suited to addressing such questions.

However, as we shall see, in trying to construct a ``topological'' object
realizing a given ``algebraic'' invariant, we will need to apply
the constructions provided in this paper to objects in both types of
category, which is why we set up our machinery in a form suitable for
both contexts.
\end{remark}

\subsection{A space and its \Pa}\label{sspace}
We now give an example of a \Pa\ which will be used later to
illustrate the general theory.

For \w[,]{k\geq n} let \w{\PAnk} denote the category of
$k$--truncated and $(n{-}1)$--connected \Pa s  $\Lambda$, with
\w{\Lambda_{i}=0}
for \w{i<n} or \w[.]{i>k} Note that in the stable range -- that
is, if \w{k<2n} -- this is an abelian category.
By restricting attention to $(n{-}1)$--connected spaces, and
truncating higher homotopy groups, we may (and shall) assume
that \w{\tru{k}\pis\bX} takes values in \w[.]{\PAnk} More formally, we
may work in the context of \fullref{stcat}(c)--(d) below.

  From now on, we take \w{n\geq 4} with \w[,]{k:=n+2} and let
\w{\gS{r} := \pis\bS{r}} and
\w{\gS{r}_{x}:=\tru{n+2}\gS{r}} denote the free monogenic algebra (in
\w{\PAlg} or \w[)]{\PAn} on a generator $x$ in degree $r$\vsm .

For \w[,]{n\geq 4} let
\w[.]{\bX:=\bS{n}\cup_{2}\be{n+1}=\Sigma^{n-1}\RR P^{2}} Then
$$
\pi_{i}\bX\cong
\begin{cases}(\ZZ/2)\lra{\alpha} & \text{for } i=n\\
(\ZZ/2)\lra{\alpha\circ\eta} & \text{for } i=n+1\\
(\ZZ/4)\lra{\beta} & \text{for } i=n+2\\
(\ZZ/2)\lra{\alpha\circ\nu}\oplus(\ZZ/2)\lra{\beta\circ\eta} &
\text{for } i=n+3
\end{cases}
$$
\noindent with \w[.]{2\beta=\alpha\circ\eta^{2}}
Note that the inclusion \w[,]{\varphi\co \tru{n+2}\pis\bX\to\gS{n-1}}
defined by \w{\varphi(\alpha)=\eta} (and \w[,]{\varphi(\beta)=6\nu}
necessarily), is a morphism of $(n{+}2)$--truncated \Pa s (in fact,
even of
$(n{+}3)$--truncated \Pa s, if \w[).]{n\geq 5}

\begin{remark}\label{rreal}
There is another non-trivial map of (truncated) \Pa s
\w[,]{\psi\co \pis\bX\to\gS{n-1}} defined by \w{\varphi(\alpha)=0}
and \w[.]{\varphi(\beta)=\eta^{3}=12\nu} This is induced by a map
of spaces -- namely, the composite of the pinch map
\w{p\co \bX=\bS{n}\cup_{2}\be{n+1}\to\bS{n+1}} with
\w[.]{\eta^{2}\co \bS{n+1}\to\bS{n-1}}
\end{remark}

%
%
\section{Resolution model categories}
\label{crmc}

In order to study the realization questions mentioned in the
Introduction, we need a suitable \emph{resolution model category}
structure on the associated simplicial model category \w[,]{\sC}
originally defined by Dwyer, Kan and Stover in \cite{DKStoE}, and
later extended by Bousfield in \cite{BouC}. A variant, called
a \emph{spiral} model category, is defined by Baues in
\cite[Section~D.2]{BauCF}. We begin with some definitions:

\begin{defn}\label{dmatch}
Let \w{(-)\otimes(-)\co \sC \times s\Sets\to\sC} be the action
of simplicial sets on the simplicial category $\sC$ (see
\fullref{dsimpst} or Quillen \cite[Section~II.1]{QuiH}).

For any finite simplicial set $K$, the \emph{matching functor}
\w{M_{K}\co \sC\to\calC} is characterized as a right adjoint by the relation
\[
\Hom_{\sC}(\cons{Z} \otimes K, \Xd) \cong
\Hom_{\calC}(Z, M_{K}\Xd).
\]
\noindent In particular,
\w[.]{M_{n}\Xd\DEF M_{\dD{n}}\Xd\DEF\lim_{[n]\to[k]}X_{k}}
Dually, the \emph{latching functor} \w{L_{n}\co \sC\to\calC} is defined by
\[
L_{n}\Xd:=\colim_{[k]\to[n]}X_{k}.
\]
Similarly, we may characterize \w{C_{K}\co \sC\to\calC} by means of a
right adjunction
\[
\Hom_{\sC}(\cons{Z} \wedge K, \Xd) \cong \Hom_{\calC}(Z, C_{K}\Xd),
\]
\noindent where \w{\Yd \wedge K} is the pushout in \w{\sC}
\mydiagram{
{\Yd \otimes \ast} \ar[r] \ar[d]
       & {\ast} \ar[d] \\
{(\Yd \otimes \ast) \otimes K} \ar[r] & {\Yd \wedge K}.
}

In particular, \w{C_{n}\Xd:= C_{M}\Xd} for
\w{M:=\Delt{n}/\Lambda^{0}{[n]}} and \w{Z_{n}\Xd:= C_{\bS{n}}\Xd}
(see \fullref{dsimpst}).
\end{defn}

\begin{remark}\label{rles}
There is a natural sequence
\[
Z_{n+1}\Xd\stackrel{i_{n+1}}{\to} C_{n+1}\Xd
\stackrel{d_{0}}{\to}
Z_{n}\Xd \stackrel{i_{n}}{\to} C_{n}\Xd,
\]
\noindent where the composite \w{i_{n}d_{0}} is induced by the map
$$\delta_{0}\co  \Delt{n}/\Lambda^{0}{[n]}
\to\Delt{n+1}/\Lambda^{0}{[n+1]}.$$
\end{remark}

Recall that we assume $\calC$ to be a right proper cofibrantly
generated pointed model category, and $\calA$ a set of models (that is,
cofibrant homotopy cogroup objects) in $\calC$.

\begin{defn}\label{daproj}
A map \w{p\co  X \to Y} in \w{\hC} is called \emph{$\calA$--epic}
if \w{p_{\ast}\co  [A,X]_{\calC} \to [A,Y]_{\calC}} is surjective for
each \w[.]{A\in \calA} An object \w{W \in \hC} is called
\emph{$\calA$--projective} if \w{p_{\ast}\co [W,X]_{\calC}\to[W,Y]_{\calC}}
is surjective for each $\calA$--epic map \w{p\co X\to Y} in \w[.]{\hC}
Finally, an object (respectively, map) of $\calC$ is called
\emph{$\calA$--projective} (respectively, \emph{$\calA$--epic}) if it is
so in \w[.]{\hC}
\end{defn}

\begin{defn}\label{dreedy}
\begin{enumerate}
\renewcommand{\labelenumi}{(\alph{enumi})}
\item a map \w{f\co \Xd\to\Yd} in \w{\sC} is a \emph{Reedy fibration} if
   the induced map \w{X_{n}\to Y_{n}\times_{M_{n}\Yd}M_{n}\Xd} is a
   fibration in $\calC$ for all \w[;]{n\geq 0}
\item a map $g$ in $\calC$ is an \emph{$\calA$--projective cofibration}
if $g$ is a cofibration in $\calC$, and has the left lifting property
with respect to the class of fibrations in $\calC$ which are, in
addition, $\calA$--epic.
\end{enumerate}
\end{defn}

\subsection{The resolution model category}\label{srmc}
Given $\calC$ and $\calA$ as above, a map \w{f\co \Xd\to\Yd} in \w{\sC} is

\begin{enumerate}
\renewcommand{\labelenumi}{(\alph{enumi})}
\item an \emph{$\calA$--weak equivalence} if
     \w{f_{\ast}\co [A,\Xd]_{\calC}\to[A,\Yd]_{\calC}} is a weak equivalence
     of
     simplicial groups for all \w[;]{A\in\calA}
\item an \emph{$\calA$--fibration} if $f$ is a Reedy fibration and
     \w{f_{\ast}\co [A,\Xd]_{\calC}\to[A,\Yd]_{\calC}} is a fibration of
     simplicial groups for all \w[;]{A \in \calA}
\item an \emph{$\calA$--cofibration} if the induced map
     \w{X_{n}\amalg_{L_{n}\Xd}L_{n}\Yd \to Y_{n}} (\fullref{dmatch}) is an
     $\calA$--projective cofibration in $\calC$ for all \w[.]{n\geq 0}
\end{enumerate}

\begin{thm}\label{tone}
If $\calC$ is a pointed right proper simplicial model category with a
set of models $\calA$, then \w[,]{\sC} with the $\calA$--weak equivalences,
$\calA$--fibrations, and $\calA$--cofibrations, and the external
simplicial category structure (\fullref{dsimpst} and Quillen
\cite[Section~II.1]{QuiH}),
is a right proper simplicial model category, called the
\emph{$\calA$--resolution model category}, and denoted by \w[.]{\CA}
\end{thm}

\begin{proof}
See Jardine \cite[Theorem~2.2]{JardB}.
\end{proof}

\begin{example}\label{ermc}
If \w{\calC =\calTs} and \w[,]{\calA:=\{\bS{n}\}_{n=1}^{\infty}}
(generated by \w[),]{\bS{1}} the resulting $\calA$--resolution model
category structure on the category \w{s\calTs} of pointed simplicial
spaces is the original ``$E^{2}$--model category'' of Dwyer, Kan and
Stover \cite{DKStoE}.
\end{example}

In constructing cofibrant replacements for objects in an
$\calA$--resolution
model category, we shall have occasion to use the following:

\begin{defn}\label{dcw}
A \emph{CW complex} is an object \w{\Xd \in \CA} such that
\begin{itemize}
\item For each \w[,]{n\geq 0} \w{X_{n} \cong \bar X_{n}\amalg L_{n}\Xd}
	for some \w[;]{\bar X_{n}\in\Obj\PiA}
\item \w{d_{i}|_{\bar X_{n}} = \ast} for all \w[.]{i\geq 1}
\end{itemize}

The \emph{attaching map} \w{d_{0}|_{\wwbar{X}_{n}}\co \wwbar{X}_{n}\to
L_{n-1}\Xd}
is denoted by \w[.]{\bar d_{0}} The collection
\w{\{\wwbar{X}_{n}\}_{n=0}^{\infty}}
is called a \emph{CW basis} for \w[.]{\Xd}
It is straightforward to check that a CW complex in \w{\CA} is
$\calA$--cofibrant.
\end{defn}

\begin{defn}\label{dnat}
The $n$th \emph{natural homotopy group} of \w{\Xd\in \sC} with
coefficients in \w{A\in \calA} is defined to be
\w{\pinat_{n}(\Xd,A)\DEF\pi_{0}\map_{\sC}(A \hotimes\bS{n},\Yd)}
(cf \fullref{dsimpst}),  where \w{\Xd\to\Yd} is a Reedy fibrant
replacement of \w[.]{\Xd} It can be equivalently defined by the exact
sequence
$$
[A,C_{n+1}\Yd]_{\calC}\stackrel{(d_{0})_{\ast}}{\to}
      [A,Z_{n}\Yd]_{\calC}\to\pinat_{n}(\Xd,A) \to 0.
$$
\noindent (see May \cite[17.3]{MayS}). Denote the $\calA$--graded group
\w{(\pinat_{n}(\Xd,A))_{A \in\calA}} by
$$
\pinat_{n}(\Xd,\calA)=\pin{n}{\Xd}\,.
$$
\end{defn}

\begin{remark}
Since \w{A\in \calC} is a homotopy cogroup object, whenever \w{\Xd\in\sC}
is Reedy fibrant we may identify  \w{[A,C_{n}\Xd]_{\calC}} with
\w{C_{n}[A,\Xd]_{\calC}} (the
$n$--chains group (\fullref{dmatch}) for the simplicial group
\w[).]{[A,\Xd]_{\calC}}
\end{remark}

\begin{defn}\label{dpisx}
By applying the functors \w{[A,-]_{\calC}} for \w{A\in\calA} to a
simplicial object \w[,]{\Xd\in\sC} we obtain a simplicial group
\w{[A,\Xd]_{\calC}}, since our models are homotopy cogroup objects by
assumption. This leads to another kind of homotopy group for
\w[,]{\Xd} namely
\w[.]{\pi_{n}(\Xd,A):=\pi_{n}[A,\Xd]_{\calC}}
Write \w{\pi_{n}\piA\Xd} for the $\calA$--graded group
\w[.]{(\pi_{n}(\Xd,A))_{A\in\calA}}
\end{defn}

As shown by Dwyer, Kan and Stover \cite[8.1]{DKStoB} and, more generally,
by Goerss and Hopkins \cite[3.4]{GHopkR},
the two types of $\calA$--graded homotopy groups are related by a
\emph{spiral exact sequence}
%
\begin{multline}
\label{eqfour}
\cdots\to\Omega\pinat_{n-1}(\Xd, A)\xra{s_{n}}
\pinat_{n}(\Xd, A)\xra{h_{n}}\pi_{n}\piA\Xd\xra{\partial_{n}} \\
\Omega\pinat_{n-2}(\Xd, A)\to\cdots\to \pinat_{1}(\Xd, A)\to
\pi_{1}\piA\Xd
\end{multline}
\noindent where \w[,]{\Omega\pinat_{n}(\Xd,A)\DEF\pinat_{n}(\Xd,\csusp A)}
for \w{\csusp A} the suspension of $A$ in $\calC$.

%
%
\begin{prop}[cf Blanc, Dwyer and Goerss {\cite[Proposition~7.13]{BDGoerR}}]\label{pfive}
For any simplicial object \w[,]{\Xd\in\CA} there are natural actions
of \w{\pinat_{0}(\Xd,\calA)\cong \pi_{0}\piA\Xd} on
\w{\pinat_{n}(\Xd,\calA)}
and \w[,]{\Omega\pinat_{n}(\Xd,\calA)} making the spiral exact sequence
\eqref{eqfour} a long exact sequence of modules over
\w[.]{\pinat_{0}(\Xd,\calA)}
\end{prop}

\begin{proof}
Because \w{\bS{n}=\Delt{n}/\partial\Delt{n}} has two non-degenerate
simplices, if we set
$$\w[,]{\widehat{A\otimes\bS{n}}:=
(A\hotimes\,\Delt{n})/(A\hotimes\,\partial\Delt{n})}$$
the map of simplicial sets \w{\bS{n}\to\Delt{0}} has a
section, which induces
$$
\xymatrix{
\widehat{A\otimes\bS{n}}\ar[r]^>>>>>{i} &
A\hotimes\bS{n} \ar[r]_<<<<{p} & \ar@/_1pc/_{s}[l] A\hotimes\Delt{0},
}
$$
and thus a natural splitting
$$
\xymatrix{
\pinat_{n}(\Xd,A) \ar[r]_{p_{\#}} & \ar@/_1pc/_{s_{\#}}[l]
\pinat_{0}(\Xd,A)}
$$
\noindent for each \w{\Xd\in\sC} and \w[.]{A\in\calA}
Using the usual homotopy cogroup structure on \w{\bS{n}}
(over \w[),]{\Delt{0}} we see that \w{\pin{n}{\Xd}} is actually a
group object over \w[.]{\pin{0}{\Xd}}  Furthermore, it is abelian because
of the underlying group structure coming from the fact that each
\w{A\in\calA} is a homotopy cogroup object itself (compare
Whitehead \cite[III, Theorem~5.21]{GWhE}).
\end{proof}

\begin{remark}
\w{\Ker(p_{\#})\cong[\widehat{A\otimes\bS{n}},\Xd]}
is actually the more traditional $n$th homotopy group of \w{\Xd}
(over the base-point component).
\end{remark}

\subsection{Algebraic categories}\label{sacat}
It will be helpful to include the following ``algebraic'' examples (cf
\fullref{rtwocat}) among our candidates for $\calC$:

\begin{enumerate}\renewcommand{\labelenumi}{(\alph{enumi})}
%
\item Let \w[.]{\calC=\PAAlg$ and $\calB=\{\piA(A) \}_{A \in \calA}}  Then
    $\calC$ has the \emph{trivial model category} structure, where only
    isomorphisms are weak equivalences and all maps are both
    cofibrations and fibrations (notice this implies the suspension
    functor \w{\csusp} is the constant functor on $\ast$).
    Recall that the objects of the form
    \w{\calA(A,?)} constitute a strong generating set for
    \w{\gr_{\calA}\Sets} by
    the Yoneda lemma, and \w{F\calA(A,?)=\piA(A)} for the free functor $F$
    defined in \fullref{scpia}.	Hence, the resolution model category
    structure
    on \w{s\PAAlg} with this $\calB$ is identical to the usual model
    category structure on \w{\sC} inherited from the category of
    simplicial
    ($\calA$--graded) groups\vsm.

\item More generally, let \w{\calC=\Alg{\Theta}} be any
    \emph{FP-sketchable} variety of (graded) universal algebras,
    corepresented by an FP-theory $\Theta$ (cf Ad\'amek and Rosick\'y
    \cite[Section~1]{ARosiL} or
    Blanc and Peschke \cite[Section~1]{BPescF}): for example, the categories
    of \APa s
    (corepresented by \w[),]{\Theta=\PiA\op} Lie algebras,
    graded commutative algebras, and so on. We assume that $\Theta$ is a
    $\fG$--theory as in \cite[Section~2]{BPescF}, so that each
    \Ta\ has an underlying (graded) group structure. In this case
    we can endow $\calC$ with the trivial model category structure,
    take $\calA$ to be the set of all monogenic free \Ta s, and obtain the
    usual model category structure on \w{\sC} (cf Quillen
    \cite[Section~II.4]{QuiH}).

\item As an application of example (b) above, if \w{\calC=\Gp} and
    \w[,]{\calA=\{\ZZ\}} then \w{\CA} (where \w[)]{\sC=\calG}
    also provides a resolution model category for the homotopy theory
    of pointed connected topological spaces (cf \cite[Section~II.3]{QuiH}).
\end{enumerate}

\begin{remark}\label{rsgp}
For many purposes it is more convenient to work with $\calG$
than with \w[.]{\Ts}  When we do so, we use the simplicial group
spheres \w{\mS^{n} =F\bS{n-1}\in\calG}	for \w{n\geq 1}
(and \w[)]{\mS^{0}=G\bS{0}} as our models $\calA$. (For definitions
of the various loop group constructions on simplicial sets see,
for example, Goerss and Jardine \cite[V.6]{GJarS}.) Note that $\DD$--diagrams of
simplicial spaces are then replaced by $\DD$--diagrams of bisimplicial
groups, which are just (more complicated) diagrams of
groups, so that many constructions may be performed entrywise in
\w[.]{\Gp}
\end{remark}

\subsection{Topological categories}\label{stcat}
It is also useful to include a number of variants of the usual
category of pointed topological spaces\vsm :

\begin{enumerate}\renewcommand{\labelenumi}{(\alph{enumi})}
%
\item If \w{\calC=\calTs} in the rational model structure and
    \w{\calA:= \{\bS{n}_{\QQ}\}_{n=2}^{\infty}} (generated by
    \w[)]{\bS{2}_{\QQ}} or \w{\calC=\calTs} in the $p$--local model
    structure and
    \w[,]{\{\bS{n}_{(p)}\}_{n=2}^{\infty}} then we have resolution model
    structures on \w{s\calTs} for rational or $p$--local simply-connected
    homotopy theory\vsm.
%
\item If \w{\calC=\Spec} is an appropriate category of spectra (cf
    Mandell, May, Schwede and Shipley \cite{MMSShipM}), and
    \w{\calA:=\{\bS{n}\}_{n=-\infty}^{\infty}}
    are all sphere spectra, we have a resolution model category
    structure on \w{s\Spec} for simplicial spectra (see
    Goerss and Hopkins \cite{GHopkR,GHopkM,GHopkM2} for the details on
    this and other categories of structured ring spectra)\vsm.
%
\item Take $\calC$ to be one of the model categories for $n$--types,
    such as the $n$--cat groups of Loday \cite{LodaS} or the crossed $n$--cubes
    of Ellis and Steiner \cite{ESteH} and \w[,]{\calA:=
    \{\bS{k}\}_{k=1}^{n}} which gives
    a resolution model category structure on \w{\sC} for $n$--types of
    spaces. An alternative is to use the (left) Bousfield localization
    model
    category structure on pointed spaces (see Hirschhorn
    \cite[Sections~2.1 and~3.3]{HirM}) for
    the map \w{\ast\to\bS{n+1}} (see Dror Farjoun
    \cite[Section~1.E.1]{DroC}).
%
\item Take \w{\calC = \calTs} and \w{\calA = \{\bS{n}\}_{n=k}^{\infty}}
   (generated by \w[);]{\bS{k}}  then we have the resolution model
   structure on \w{s\calTs} for the homotopy theory of ``$(k{-}1)$--connected
   types'' for spaces -- that is, the right Bousfield localization
   model of \cite[Section~3.3]{HirM}  (see \cite[Section~2.D.2.6]{DroC}).
\end{enumerate}

\subsection{Diagram categories}\label{sdcat}
The motivating type of example for this paper was the category
\w{\calTs^{\DD}} of \emph{$\DD$--diagrams of spaces}, where $\DD$ is
a small category.

Recall that for any object \w{X\in\calC} and \w[,]{d\in\Obj\DD} the
free $\DD$--diagram \w{F(X,d)} is defined by setting the $e$--entry
equal to \w[,]{F(X,d)_{e}:=\coprod_{\Hom_{\DD}(d,e)}X}
with maps induced by the identity on each factor.
Then for any collection of models $\calA$ for $\calC$, the
\emph{induced collection of models} $\calB$ for \w{\calC^{\DD}}
consists of all free $\DD$--diagrams \w{F(A,d)} for \w{d\in\Obj\DD}
and \w[.]{A\in\calA}

Note that the model category structure on \w{s\calTs^{\DD}} given by
\fullref{tone} using $\calB$ is identical to the structure induced
from that on \w{s\calTs} associated to $\calA$ (and \fullref{tone}) as in
\cite[Section~11.6]{HirM}.  Furthermore, the category
\w{\PAAlg} is equivalent to the category of $\DD$--diagrams of
(ordinary) \Pa s in these cases.

\begin{notn}\label{nnca}
For any \w[,]{n\in\NN} let \w{[n]} denote the category
with \w{n+1} objects \w{0,1,\ldots,n} and $n$ composable maps between
them.
For example, \w{\DD=[1]} has two objects and a single non-identity
morphism \w[.]{0\to 1}
\end{notn}

\begin{examples}\label{edcat}
\begin{enumerate}\renewcommand{\labelenumi}{(\alph{enumi})}
%
\item If \w{\calC=\calTs} and \w[,]{\DD=[1]} then
\w{\calTs^{\DD}} is the category of \emph{maps of spaces}, and
for any space $X$, the free object
\w[,]{F(X,0)=X\smash{\stackrel{\Id}{\to}}X}
while \w[.]{F(X,1)=\ast\to X}  Hence in this case
\w{\calA:=\bigl\{\ast\to\bS{n},\bS{n}\smash{\stackrel{\Id}{\to}}\bS{n}\bigr\}_{n=1}^{\infty}}
-- that is, $\calA$ is generated by the pair consisting of
\w{\ast\to\bS{1}} and \w{\bS{1}\smash{\stackrel{\Id}{\to}}\bS{1}} -- and
\w{\PAAlg} is the category of morphisms between \Pa s\vsm.
%
\item Suppose \w{\calC=\calTs} and \w{\DD=[2]} (with a single
    composable pair of nonidentity maps, denoted \w[).]{0\to 1\to 2}
    Then for any space $X$,
    \w[,]{F(X,0)=X\smash{\stackrel{\Id}{\to}}X\smash{\stackrel{\Id}{\to}}X}
    \w[,]{F(X,1)=\ast\to X\smash{\stackrel{\Id}{\to}}X} and
\w[.]{F(X,2)=\ast\to\ast\to X} Thus $\calA$ is generated by:
$$
\ast\to\ast\to\bS{1},\hsp \ast\to\bS{1}\stackrel{\Id}{\to}\bS{1},\hsp
\text{and }\bS{1}\stackrel{\Id}{\to}\bS{1}\stackrel{\Id}{\to}\bS{1}.
$$
\noindent while \w{\PAAlg} is the category of composable pairs of maps
between \Pa s.
\end{enumerate}
\end{examples}

%
%
\section{$E^{2}$--model categories}
\label{cemc}

There are a number of familiar constructions for topological spaces
which relate to Postnikov towers and are useful to have in a
resolution model category \w[,]{\CA} although they need not exist
in general.
We shall show, however, that these constructions are available in all
of the examples we wish to consider.

\begin{defn}\label{dfpt}
A \emph{Postnikov tower} functor applied to an object \w{\Xd} in a
resolution model category \w{\CA} is a functorial commuting diagram
\begin{equation}{\label{epost}}
\vcenter{\xymatrix@C=16pt{
\Xd \ar[dr]^>>>>{r\q{n+1}} \ar@/^1pc/[drr]^>>>>{r\q{n}}
\ar@/^2pc/[drrr]^>>>>{r\q{n-1}} & & & & \\
\cdots \ar[r] & ~~P_{n+1}\Xd~~~~ \ar[r]^{p\q{n+1}} &
~~P_{n}\Xd~~~~ \ar[r]^{p\q{n}} &
~~P_{n-1}\Xd~~~~ \ar[r]^{p\q{n-1}} & ~~~~\cdots P_{0}\Xd
}}
\end{equation}
\noindent of $\calA$--fibrations \w{p\q{n}} and maps \w{r\q{n}}
which induce isomorphisms
\[
\pinat_{k}(P_{n}\Xd;\calA)\cong
\begin{cases}
	       \pinat_{k}(\Xd;\calA) &0 \leq k \leq n;\\
	       0		     &\text{otherwise}.
\end{cases}
\]
\end{defn}

\begin{defn}\label{drem}
If \w{\CA} is a resolution model category, a \emph{classifying object}
\w{\BL=\BCL{\sC}} for a \APa\  $\Lambda$ is any fibrant \w{\Bd\in\sC}
such that \w{\Bd\simeq P_{0}\Bd} and \w[.]{\pin{0}{\Bd}\cong\Lambda}
\end{defn}

\begin{defn}\label{dem}
Given an abelian \APa\ $M$ and an integer \w[,]{n\geq 1} an
\emph{$n$--dimensional $M$--Eilenberg--Mac\,Lane object}
\w{\EM{}{M}{n}=\EC{\sC}{M}{n}} is any fibrant \w{\Ed\in\sC} such
that
\w{\pin{n}{\Ed}\cong M} and \w{\pin{k}{\Ed}=0} for \w[.]{k\neq n}
\end{defn}

\begin{defn}\label{deaem}
Given a \APa\ $\Lambda$, a module $M$ over $\Lambda$, and an integer
\w[,]{n\geq 1} an \emph{$n$--dimensional extended $M$--Eilenberg--Mac\,Lane
object}
\w{\EL{M}{n}=\ECL{\sC}{M}{n}} is any fibrant homotopy abelian group
object \w{\Ed\in\sC/\BL} satisfying
%
%
\begin{equation}\label{eqone}
\pin{k}{\Ed}\cong \begin{cases}
	     \Lambda & \text{for }k=0,\\
	     M \text{ (as a module over $\Lambda$) }& \text{for } k=n,\\
	     0	   & \text{otherwise}. \end{cases}
\end{equation}
\end{defn}

\begin{remark}\label{rsec}
The fact that \w{\Ed=\EL{M}{n}} is a homotopy abelian group
object in \w{\sC/\BL} implies that
\w{[\BL,\Ed]_{\sC/\BL}} has a natural abelian group
  structure, so in particular a unit element.
Thus \w{\Ed} comes equipped with a designated homotopy section $s$ for
\w[.]{r\q{0}\co \Ed\to P_{0}\Ed\simeq\BL}
\end{remark}

From the spiral exact sequence \eqref{eqfour} we readily calculate
%
\begin{equation}\label{eten}
\pi_{k}\piA\EL{M}{n}\cong \begin{cases}
		       \Lambda & \text{for }k=0,\\
		       \Omega\Lambda & \text{for }k=2,\\
		       M & \text{for }k=n,\\
		       \Omega M & \text{for }k=n+2,\\
		       0     & \text{otherwise}, \end{cases}
\end{equation}
\noindent with the obvious variant when \w{n=2} (that is,
\w[).]{\pi_2\piA\EL{M}{2} \cong \Omega\Lambda\times M}

\begin{remark}\label{rloop}
Note that if we apply the loop functor in the category \w{\sC/\BL}
to \w{\EL{M}{n}} -- that is, take the pullback of
\w{\BL\leftarrow \EL{M}{n}\to\BL} (cf Quillen \cite[Section~I.2]{QuiH}) --
we obtain \w[.]{\EL{M}{n-1}}
\end{remark}

\begin{defn}\label{dki}
Given a Postnikov tower functor as in \fullref{dfpt}, an $n$th
\emph{$k$--invariant square} (with respect to $\calA$)
is a functor that assigns to each \w{\Xd\in\sC} a homotopy pull-back
square
%
%
\begin{equation}{\label{eqfive}}
\vcenter{\xymatrix@C=35pt{
\ar @{} [dr] |<<<<<<{\framebox{\scriptsize{hPB}}}
P_{n+1}\Xd \ar[r]^{p\q{n+1}} \ar[d] &
       {P_{n}\Xd}    \ar[d]^{k_{n}}\\
\BL \ar[r]_<<<<<<<{s} & {\EL{M}{n+2}}
}}
\end{equation}
\noindent for \w{\Lambda:=\pin{0}{\Xd}} and \w[.]{M:=\pin{n+1}{\Xd}}
The map \w{k_{n}\co P_{n}\Xd\to\EL{M}{n+2}} is called the $n$th
$k$--\emph{invariant} for \w[.]{\Xd}
\end{defn}

\begin{defn}\label{demc}
A resolution model category \w{\CA} as in \fullref{srmc} is called an
\emph{$E^{2}$--model category} if:

\begin{enumerate}
\renewcommand{\labelenumi}{Ax~\arabic{enumi}}
%
\item\label{ax1} $\sC$ has functorial Postnikov towers\vsm.
%
\item\label{ax2} For every \APa\ $\Lambda$ and $\Lambda$--module $M$  the
    classifying object \w{\BL} and the $n$--dimensional extended
    Eilenberg--Mac\,Lane object \w{\EL{M}{n}} exist, for each
    \w[.]{n\geq 1} In addition we assume the latter determines a functor
$$
\EL{-}{n} \co   \RM{\Lambda} \to\ab(\ho(s\calC)),
$$
both constructions are functorial in $\Lambda$, and are unique up to
homotopy.
\item\label{ax3} $\sC$ has $k$--invariant squares (with respect to $\calA$) for
    each \w[.]{n\geq 0\vsm}
\item\label{ax4} There is a functor \w{J\co \sC\to\calC} such that, for
\w{\Lambda\in\PAAlg}
and \w[,]{\Xd\in\sC} if
\w{\piA\Xd\smash{\stackrel{\hbox{\footnotesize$\sim$}}{\to}}B_{s\PAAlg}\Lambda}
is a weak equivalence in \w[,]{s\PAAlg} then there is an isomorphism
%
\begin{equation}\label{eqsix}
[A,J\Xd]_{\calC} \stackrel{\cong}{\rightarrow}\Hom_{\PAAlg}(\piA
A,\Lambda),
\end{equation}
\noindent natural in $\Lambda$ and \w[.]{A\in\calA}
\end{enumerate}
\end{defn}

\begin{remarks}\label{rrss}\quad

\begin{enumerate}
\renewcommand{\labelenumi}{$\bullet$}
\item Ax~\ref{ax1}--Ax~\ref{ax3} imply that \w{\CA} is a \emph{spherical model category}
in the
sense of Blanc \cite[Section~2]{BlaCH}, and so in particular is
\emph{stratified}
in the sense of Spali\'nski \cite{SpalSM}. These axioms are also satisfied, for
example,
by the category \w[,]{\calTs} which is not itself a resolution model
category
(but see \fullref{rsgp}).
\item We may assume that our extended Eilenberg--Mac\,Lane objects
are \emph{strict} abelian group objects in \w[,]{\sC/\BL} by
functoriality,
since the group structure morphisms for a $\Lambda$--module $M$ are
maps of modules.
\item Not all resolution model categories have the additional
  structure of a spherical model category (see \fullref{retwo})\vsm .
\item The point of Ax~\ref{ax4} is that any \w{\Xd \in \sC/\BL} with
\w{\piA\Xd\simeq B_{s\PAAlg}\Lambda} in \w{s\PAAlg}
yields a realization \w{J\Xd} for $\Lambda$ (see \fullref{tfour}).
See Chach\'olski, Dwyer and Intermont \cite{CDIComp} for a way to
geometrically handle cases where Ax~\ref{ax4} does not hold.
\item The statement of Ax~\ref{ax4} may appear somewhat convoluted, because it
  is intended to apply to two rather different contexts: see
  Theorems~\ref{ttwo} and~\ref{tfive} below.
 \fullref{ttwo} deals with the case of universal algebras (hence
  the special case of \APa s), while \fullref{tthree} treats the
  general extension to diagram categories, thereby reducing  our
  motivating example of diagrams of spaces to a consequence of \fullref{tfive}, which deals with \w{s\calTs} with several standard
  model structures on \w[.]{\calTs}
\end{enumerate}
\end{remarks}

%
%
\begin{thm}\label{ttwo}
Let \w{\calC=\Alg{\Theta}} be an FP-sketchable variety of (graded)
universal
algebras, corepresented  by a $\fG$--theory $\Theta$,  with trivial
model category structure,  and let $\calA$ consist of monogenic free
\Ta s, as in \fullref{ermc}(f). Then \w{\CA} is an $E^{2}$--model
category.
\end{thm}

\begin{proof}
We use the constructions described by Blanc, Dwyer and Goerss
\cite{BDGoerR} for the case
$\calC=\PAlg$\vsm:

\noindent \textbf{For Ax~\ref{ax1}}\qua Follow Dwyer and Kan
\cite[Section~1.2]{DKanO}:

Given \w{\Yd\in\sC} and \w[,]{n\geq 0} first define \w{\Yd\q{n}\in\sC} by
\[
Y^{(n)}_{k}=\begin{cases}
	Y_{k}		&0 \leq k \leq n+1;\\
	M_{k}(\Yd\q{n}) &n+2 \leq k,
\end{cases}
\]
\noindent with simplicial maps determined from	\w{\tru{n+1}\Yd}
and \w[,]{\delta_{k} \co  M_{k}(\Yd^{(n)}) \to Y_{k}^{(n)}} along with the
obvious maps \w{p\q{n}\co \Yd\q{n}\to\Yd\q{n-1}} and
\w[.]{r\q{n}\co \Yd\to\Yd\q{n}}

The Postnikov tower for \w{\Xd\in\sC} is then defined by setting
$P_{n}\Xd:=\Yd^{(n)}$, where \w{\Xd\to\Yd} is a (functorial)
$\calA$--fibrant replacement in \vsm\w[.]{\CA}

\noindent \textbf{For Ax~\ref{ax2}}\qua Follow \cite[Proposition~2.2]{BDGoerR},
taking \w{\BL} to be the constant simplicial object on
$\Lambda$, \w{\EM{}{M}{n}} to be the iterated Eilenberg--Mac\,Lane
construction $\wwbar{W}$ on \w{BM} (cf May \cite[Section~20]{MayS}),
and \w{\EL{M}{n}} to be the semi-direct product \w{\BL\ltimes\EM{}{M}{n}}
(\fullref{smod})\vsm.

More explicitly, let $W$ be a free \Ta\  equipped with a surjection
\w[.]{\phi\co W\to M} Define a simplicial \Ta\ \w{\Bd} by setting
\w{\sk{n-1}\Bd:=\sk{n-1}\BL} and \w[,]{E_{n}\simeq W\amalg\BL_{n}}
with \w[.]{W\subseteq Z_{n}\Bd} A straightforward calculation shows
\w[,]{C_{n}\BL=Z_{n-1}\BL=0}  so \w{Z_{n}\Bd=C_{n}\Bd} is the cokernel
\w{\BL_{n}\ltimes W} of \w[.]{\BL_{n}\to E_{n}=W\amalg\BL_{n}}
Note that \w{\BL_{0}} embeds in \w{\BL_{n}} as a free retract by
\w[,]{s_{n-1}\dotsb s_{0}} so \w{\BL_{n}\cong\BL_{0}\amalg L'} for
some \Ta\ \w[,]{L'} where \w{L'\ltimes W} is a \Ta\ ideal in
\w[,]{Z_{n}\Bd} with quotient \Ta\
\w[.]{Z_{n}\Bd/(L'\ltimes W)\cong K_{0}\ltimes W} This is by definition
the free $\BL_{0}$--algebra generated by $W$, and thus \w{\phi\co W\to M}
extends to a map of $\BL_{0}$--algebras
\w[;]{\hat{\phi}\co \BL_{0}\ltimes W\to M} precomposing with the projection
\w{Z_{n}\Bd\to \BL_{0}\ltimes W} defines \w[.]{\tilde{\phi}\co Z_{n}\Bd\to M}

Let \w{\bar d_{0}\co \wwbar{B}_{n+1}\to B_{n}\Bd := \Ker\tilde{\phi}} be a
surjection from a free \Ta, let
\w[,]{\BL_{n+1} := \wwbar{B}_{n+1}\amalg L_{n+1}\Bd} and let
\w[.]{\Bd := P_{n}\sk{n+1}\Bd} Then \w{\pi_{n}\Bd\cong M}
(as a $\Lambda$--module), and \w{\pi_{i}\Bd=0} for \w[.]{i\neq 0,n} The
section is induced by the inclusion \w[.]{\sk{n+1}\BL\hra\sk{n+1}\Bd\vsm}

\noindent \textbf{For Ax~\ref{ax3}}\qua Follow \cite[Sections~5--6]{BDGoerR}\vsm.

Given \w{\Xd\in\sC/\BL} and \w[,]{n\geq 0} take the pushout
$$
\xymatrix{P_{n+1}\Xd \ar[rr]^{p\q{n+1}} \ar[d]
\ar @{} [drr] |>>>>>>>{\framebox{\scriptsize{PO}}} & & P_{n}\Xd\ar[d]^{f}
\\
\BL \ar[rr]_{g} & & \Yd,
}
$$
\noindent and apply the functor $P_{n+2}$ to the resulting
diagram. The connectivity argument of \cite[Lemma~5.11]{BDGoerR}
applies here, too, so the result is actually a homotopy pull-back square,
\w{P_{n+2}\Yd} is an extended Eilenberg--Mac\,Lane object
(with section \w[),]{P_{n+2}g} and  \w{P_{n+2}f} is the
$k$--invariant. The construction is evidently
natural, since we have natural Postnikov systems\vsm.

\noindent \textbf{For Ax~\ref{ax4}}\qua Use \w{\pi_{0}\co \sC\to\calC} as the
functor $J$.  Then the trivial model category structure on $\calC$
gives the first identity
\[
[A,J\BL]_{\calC}=\Hom_{\calC}(A,\pi_{0}\BL)\cong\pi_{0}\BL(A)
\]
and the second isomorphism comes from the fact that $A$ is monogenic
free, while
\w{\pi_{0}\BL\cong\pinat_{0}(\BL)\cong\Lambda} completes the claim\vsm.
\end{proof}

%
%
\begin{thm}\label{tthree}
Let \w{\sC_{\calA}} be an $E^{2}$--model category, $\DD$ a small
category, and
$\calB$ the induced collection of models in \w{\calC^{\DD}} (see
\fullref{sdcat}); then \w{(\sC^{\DD})_{\calB}} is an $E^{2}$--model
category.
\end{thm}

\begin{proof}
We use the induced collection of models $\calB$ (\fullref{sdcat}) to
extend the $E^{2}$--model structure to \w[.]{\sC^\DD} The
underlying simplicial model category structure on \w{\calC^\DD} is
that of Hirschhorn \cite[Section~11.6]{HirM}, with  weak equivalences
and fibrations defined objectwise; thus evaluation at \w{d\in\Obj\DD}
preserves fibrations and weak equivalences and forms part of a strong
Quillen pair, with left adjoint \w{F(-,d)} (the free diagram functor at
$d$). See \cite[11.5.26]{HirM}.

Hence, for \w[,]{A\in \calA} \w[,]{d\in \DD} and \w[,]{X \in \sC^\DD}
we have a natural isomorphism
\begin{equation}\label{evaleq}
[F(A,d),X]_{\sC^\DD} \cong [A,X(d)]_{\sC}.
\end{equation}
In particular, \w{\pi_{\calB}(-,F(A,d))} is the same as \w{\piA(-,A)}
after pre-composition with evaluation at $d$. By
the spiral exact sequence \eqref{eqfour}, the same holds for
\w[.]{\pinat_{\ast}(-,\calB)}

The axioms of \fullref{demc} can therefore be verified by
applying the various constructions of \w{\sC} at each $d$ in $\DD$,
and checking that the  requisite properties are satisfied in
\w[,]{\sC^\DD} once they hold objectwise\vsm:

\noindent \textbf{For Ax~\ref{ax1}}\qua Since \w{\sC} has functorial
Postnikov towers, \w{\sC^{\DD}} possesses such towers, with
\w[.]{(P_n\Xd)(d) = P_n(\Xd(d))\vsm}

\noindent \textbf{For Ax~\ref{ax2}}\qua Given a $\PiB$--algebra $\Lambda$
(that is, a functor \w[)]{\Lambda\co \DD\to\PAAlg} and a module $M$ over
$\Lambda$, for each \w{n\geq 1} we define the classifying object
\w{\BL} and extended $M$--Eilenberg--Mac\,Lane object \w{\EL{M}{n}}
objectwise, by applying the appropriate functors in \w{\sC} to the
diagrams $\Lambda$ and $M$. This is evidently functorial, unique up
to homotopy, and satisfies \eqref{eqone}. Note that in order for
\w{\EL{M}{n}} to be a homotopy abelian group object in
\w[,]{\sC^\DD/\BL} we must produce structure maps
\begin{equation}\label{eqtwo}
\mu\co \EL{M}{n}\times_{\BL}\EL{M}{n}\to\EL{M}{n},
\iota\co \EL{M}{n}\to\EL{M}{n}
\end{equation}
\noindent (over \w[),]{\BL}  satisfying the appropriate identities.
(The unit element is represented by the section \w[.)]{s\co \BL\to\EL{M}{n}}
However, since $M$ is itself an abelian group object in
\w[,]{\PAAlg/\Lambda}
it is equipped in turn with maps
$$
m\co M\times_{\Lambda}M\to M \text{ and }
i\co M\to M
$$
\noindent in \w[,]{\PAAlg/\Lambda} which are themselves maps of
$\Lambda$--modules, and these induce the maps of \eqref{eqtwo} by
functoriality. Note that the functors \w{\EL{-}{n}} in
\w{\sC} preserve products of modules (over $\Lambda$) because of the
homotopy uniqueness and functoriality\vsm.

\noindent \textbf{For Ax~\ref{ax3}}\qua Since Postnikov towers and extended
Eilenberg--Mac\,Lane objects, as well as fibrations and weak equivalences
are defined object-wise for \w[,]{d\in\Obj\DD}
defining $k$--invariants in \w{\sC^\DD/\BL} objectwise will give
homotopy pullback squares that are $k$--invariant squares\vsm.

\noindent \textbf{For Ax~\ref{ax4}}\qua Suppose we are given a functor
\w{J\co \sC\to\calC} with the requisite properties. Define
\w{J^\DD\co \sC^\DD \to \calC^\DD} by \w[.]{(J^\DD\Xd)(d) = J(\Xd(d))}
Let \w{\piA\Xd\smash{\stackrel{\hbox{\footnotesize$\sim$}}{\to}}B_{s(\PAAlg)^{\DD}}\Lambda} be a weak
equivalence. Now we have two natural isomorphisms
$$
[F(A,d), J^\DD(\Xd)]_{\calC^\DD} \cong [A, J(\Xd (d))]_{\calC}
$$
and
$$
[\pi_{\calB} F(A,d),\Lambda]_{(\PAAlg)^\DD}\cong[\piA A,\Lambda
(d)]_{\PAAlg}.
$$
 From Ax~\ref{ax4}, applied to \w{\piA\Xd (d)
 \smash{\stackrel{\hbox{\footnotesize$\sim$}}{\to}}
B_{s\PAAlg}\Lambda(d)}
in \w[,]{s\PAAlg} we have the natural isomorphism
$$
[A,J(\Xd (d))]_{\calC} \stackrel{\cong}{\to} [\piA A, \Lambda
(d)]_{\PAAlg}.
$$
Combining all three isomorphisms gives the required natural isomorphism
$$[F(A,d), J^\DD(\Xd)]_{\calC^\DD} \stackrel{\cong}{\to} [\pi_{\calB}
   F(A,d), \Lambda]_{(\PAAlg)^\DD}.\proved$$
\end{proof}

%
%
\begin{thm}\label{tfive}
The category \w{s\calTs} of simplicial pointed connected topological
spaces (with
the spheres \w{(\bS{n})_{n=1}^{\infty}} as models), and the four
examples of \fullref{stcat}, are all $E^{2}$--model categories.
\end{thm}

\begin{proof}
The case \w{\calC=\calTs} was treated in \cite{BDGoerR}, and all five
cases may be treated similarly\vsm:

\noindent \textbf{For Ax~\ref{ax1}}\qua As in the proof of \fullref{ttwo}\vsm.

\noindent \textbf{For Ax~\ref{ax2}}\qua Follow \cite[7.7]{BDGoerR}\vsm.

More explicitly, given \w[,]{A \in \calA} for each \w{n \geq 1}
recall \w[,]{\pinat_{n}(\Xd,\calA)\cong[{A} \hotimes \bS{n},\Xd]_{\sC}}
where
\w{{A} \hotimes \bS{n}} denotes \w{\ssusp{\cons{A}}{n}\in\sC} (see also
\fullref{dsimpst}).

For the existence of \w[,]{\BL} let \w{U,V \in\PiA} be such that
\w{\piA U \to\Lambda} is a free cover of $\Lambda$, and \w{\piA
V\to\piA U}
covers minimally the corresponding relations. For each summand $A$ in $V$,
attach a copy of \w{{A} \hotimes \bS{n}} to $U$. Applying \w{P_{0}} to the
resulting object of \w{\sC} yields a classifying object \w{\BL} as
required\vsm .

For the Eilenberg--Mac\,Lane objects, again we follow \cite[7.7]{BDGoerR}:

Let $W$ be  the model for \w{\BL} constructed as above. Let
\w{U,V\in\PiA} be such that \w{\piA V\to\piA U\to M} is
a presentation for $M$. Attach a copy of \w{{A} \hotimes \bS{n}} for each
summand $A$ of $U$ to form an object \w[,]{Z\in\sC} then attach a
copy of \w{{A} \hotimes \bS{n+1}} to $Z$ for each $A$--coproduct summand of
$V$ to form \w[.]{Z'} Applying \w{P_{n}} to \w{Z'} yields the
desired \w[.]{\EL{M}{n}} The existence of the section
\w{\sigma\co \BL\to\EL{M}{n}} follows from \cite[Proposition~4.9]{BDGoerR}\vsm.

\noindent \textbf{For Ax~\ref{ax3}}\qua Again follow \cite[Sections~5--7]{BDGoerR},
with the same construction as in the proof of Ax~\ref{ax3} for \fullref{ttwo}\vsm.

\noindent \textbf{For Ax~\ref{ax4}}\qua For the standard model of
\w[,]{\calC=\calTs} $J$ will be the realization or diagonal functor
\w{\|-\|\co \sC\to\calC} (left adjoint to the constant
functor \w[).]{c(-)_{\bullet}\co \calC\to\sC}
This extends entrywise to diagrams of simplicial spaces,
as does the natural spectral sequence of Quillen \cite{QuiS} (see also
Bousfield and Friedlander \cite[Theorem~B.5]{BFrH}), yielding an
$(\NN\times\calA)$--graded spectral sequence with
%
\begin{equation}\label{eqnine}
E^{2}_{s,\calA}=\pi_{s}(\Xd,\calA)\Rightarrow\piA\|\Xd\|.
\end{equation}
\noindent Then \eqref{eqsix} will be the edge homomorphism of this
spectral
sequence, which collapses at the $E^{2}$--term if
$\piA\Xd\simeq\piA\BL$\vsm.

\noindent We can extend this spectral sequence argument to the other
examples of \fullref{stcat} as follows:

\begin{enumerate}
\renewcommand{\labelenumi}{(\roman{enumi})}
\item For \fullref{stcat}~(a): the exactness of \w{-\otimes R} for
\w{R\subseteq\QQ} allows us to obtain a localized Quillen spectral
sequence to verify Ax~\ref{ax4} for either rational or $p$--local spaces.
\item For \fullref{stcat}~(b): the spectral sequence for the realization
   of a simplicial spectrum is analyzed by Goerss and Hopkins in
   \cite[Section~6]{GHopkR}, showing
   that Ax~\ref{ax4} is satisfied for \w{s\Spec} (as  well as for some
   structured versions of spectra). For the remaining axioms see
   \cite{GHopkM,GHopkM2}.
\item For \fullref{stcat}~(c): to verify Ax~\ref{ax4}, apply the
    Quillen spectral sequence to \w[,]{P_{n}\Xd} where \w{\Xd} is the
    usual resolution in \w[.]{s\calTs} Note that
    \w{P_{n}\|\Xd\|} is $n$--equivalent to \w{\|P_{n}\Xd\|} (as we
    can see from the differentials in the spectral sequence itself).
\item For \fullref{stcat}~(d): if \w[,]{\calA:=\{\bS{n}\}_{n=k}^{\infty}}
    we can use the usual Eilenberg--Mac\,Lane objects (noting that the
    connectivity
    assumptions are not in the simplicial direction), and again apply the
    Quillen spectral sequence to resolutions in which all spaces happen to
    be $(k{-}1)$--connected.\proved
\end{enumerate}
\end{proof}

\begin{remark}\label{retwo}
Note that not all resolution model categories are $E^{2}$--model
categories. In particular, if we replace the spheres by Moore
spaces as our models (in $\calTs$), then we have neither
Eilenberg--Mac\,Lane
objects nor Postnikov systems for the mod $p$ homotopy groups (see
Blanc \cite[Section~3.10]{BlaCH}). In addition, the realization of simplicial spaces
does not provide the expected functor  $J$ for Ax~\ref{ax4}, since the
Bousfield--Friedlander spectral sequence for a mod $p$ resolution does
not collapse (see Blanc \cite[Section~4.6]{BlaM}).
\end{remark}

\begin{notn}\label{npalgsp}
In what follows we will often have to deal with parallel constructions
of the $E^{2}$--model category structure in \w[,]{\CA} as well as in
the associated algebraic category \w[.]{s\PAAlg} In order to
distinguish between them, we shall use boldface -- \w[,]{\bP{n}\Xd}
\w[,]{\bBL:=\BCL{\sC}} \w[,]{\bE{}{M}{n}:=\EC{\sC}{M}{n}} and so on --
for the constructions in \w[,]{\sC} and tildes -- \w[,]{\tP{n}\Gd}
\w[,]{\tBL:=\BCL{s\PAAlg}} \w[,]{\tE{}{M}{n}:=\EC{s\PAAlg}{M}{n}}
etc. -- for the analogous constructions in \w[.]{s\PAAlg}

We may still use the unadorned symbols \w[,]{P_{n}\Xd} \w[,]{\BL} and
\w[,]{\EL{M}{n}} etc., when we do not need to make this distinction.
\end{notn}

%
%
\section{Cohomology theories}
\label{ccoh}

As one might expect, the Eilenberg--Mac\,Lane objects in an
$E^{2}$--model category can be used to define suitable cohomology theories:

\begin{defn}\label{deist}
Let \w{\CA} be any resolution model category. A sequence of pointed
contravariant functors \w{(D^{n}\co \ho\CA\to\RM{\ZZ})_{n=0}^{\infty}}
is called a sequence of \emph{cohomology functors} if they satisfy the
analogues of the usual Eilenberg--Steenrod axioms\vsm:
\begin{enumerate}
\renewcommand{\labelenumi}{\Roman{enumi}}
\item \w{D^{n}(\coprod_{\alpha} X_{\alpha})\cong
	       \prod_{\alpha} D^{n}X_{\alpha}}
      for any coproduct of cofibrant objects in \w[.]{\CA}
\item \w{D^{i}({A} \hotimes \bS{n})=0} for \w{i\neq n} and any
\w[;]{A\in\calA}
\item Given \w{\Nd\leftarrow \Md\xra{i}\Pd} in \w[,]{\sC} with all
    objects cofibrant and $i$ a cofibration, let
    \w{\Xd:=\Nd\amalg_{\Md}\Pd} be the pushout.  Then there is
    a natural \emph{Mayer--Vietoris} long exact sequence
%
\begin{multline}\label{eqten}
0\to D^{0}\Xd \to D^{0}\Nd\oplus D^{0}\Pd \to
D^{0}\Md\to D^{1}\Xd\\
\cdots \to D^{n}\Xd\to D^{n}\Nd\oplus D^{n}\Pd \to D^{n}\Md\to\cdots
\end{multline}
\end{enumerate}
\end{defn}

\begin{defn}\label{dcoh}
Fix a \APa\ $\Lambda$ and a $\Lambda$--module $M$. For
\w{\Xd\in\CA/\BL} and \w[,]{n\geq 1} define the \emph{$n$th
  (andr\'{e}-Quillen) cohomology group of \w{\Xd} over $\Lambda$ with
  coefficients in $M$}, denoted by \w[,]{H^{n}_{\Lambda}(\Xd;M)} to be
\[
H^{n}_{\Lambda}(\Xd;M) :=[\Xd,E_{\Lambda}(M,n)]_{\CA/\BL}.
\]
\end{defn}

We would like to know that extending \w{\piA\co \CA/\bBL\to s\PAAlg/\piA\bBL}
to a functor \w{\CA/\bBL \to s\PAAlg/\tBL} (via
$\piA\bBL\to\tP{0}\piA\bBL\simeq\tBL$) induces an isomorphism of
cohomology theories over $\Lambda$. This holds for \w{n \geq 2}
by the following generalization  of Blanc, Dwyer and Goerss
\cite[Proposition~8.7]{BDGoerR}:

%
%
\begin{prop}\label{psix}
There is a natural map \w{\zeta\co \piA\bEL{M}{n}\to\tEL{M}{n}} such that
\[
\phi_{n}(\Xd)\co [\Xd,\bEL{M}{n}]_{\CA/\bBL}\to[\piA\Xd,\tEL{M}{n}]_{s\PAAlg/\tBL},
\]
defined as the composite of the maps induced by $\zeta$ and
\w[,]{\piA\co \sC\to s\PAAlg} is an isomorphism for \w[.]{n \geq 2}
\end{prop}

\begin{proof}
The section \w{\sigma\co \bBL\to\bEL{M}{n}} (\fullref{rsec}) induces a
section \w{s\co \piA\bBL\to\tP{n}\piA\bEL{M}{n}} for the map
\w{\tilde{p}\q{n}\co \tP{n}\piA\bEL{M}{n}\to\tP{n-1}\piA\bEL{M}{n}=\piA\bBL}
(cf \fullref{dfpt}) over \w[.]{\tBL} Moreover, \w{\piA\bEL{M}{n}}
is known from \eqref{eten}. Therefore, the $(n{-}1)$st
$k$--invariant for \w{\piA\bEL{M}{n}} fits into a homotopy-commutative
diagram
$$
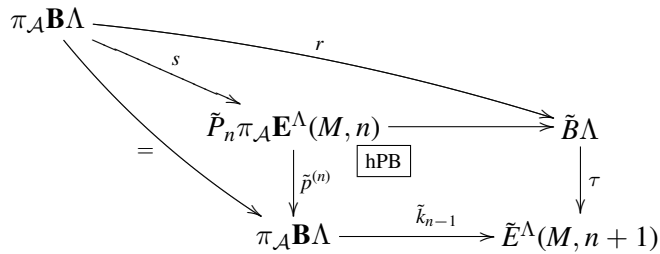
\xymatrix@C=40pt{
\piA\bBL  \ar@/_/[ddr]_{=} \ar[dr]^{s} \ar@/^/[drr]^{r}&  & \\
    & \ar @{} [dr] |<<<<<{\framebox{\scriptsize{hPB}}}
   \tP{n}\piA\bEL{M}{n} \ar[d]^{\tilde{p}\q{n}} \ar[r]	& \tBL
   \ar[d]^{\tau} \\
    & \piA\bBL \ar[r]^{\tk{n-1}}  & \tEL{M}{n+1}
}
$$
\noindent where \w{\tilde{p}\q{n}} is induced by
\w[,]{\piA(\bp{n})\co  \piA\bEL{M}{n} \to \piA\bBL} and $r$ and
the unlabelled arrow is the unique terminal map in \w[.]{s\PAAlg/\tBL}
Thus \w[,]{\tk{n-1}=\tau\circ r} yielding two consecutive homotopy
pullback squares
$$
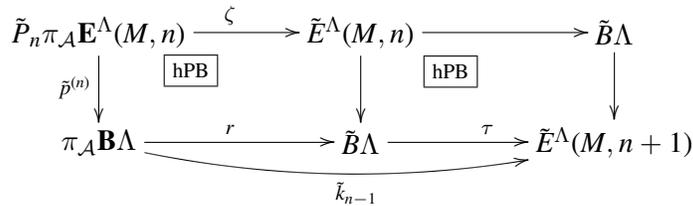
\xymatrix@C=40pt{
\ar@{}[dr] |<<<<<<{\framebox{\scriptsize{hPB}}}
\tP{n}\piA\bEL{M}{n} \ar[d]_{\tilde{p}\q{n}} \ar[r]^{\zeta} &
\ar@{}[dr] |<<<<<<{\framebox{\scriptsize{hPB}}}
\tEL{M}{n} \ar[d] \ar[r]  & \tBL \ar[d] \\
\piA\bBL \ar[r]^{r} \ar@/_1pc/[rr]_{\tk{n-1}} & \tBL \ar[r]^{\tau}  &
\tEL{M}{n+1}
}
$$
\noindent in which the required $\zeta$ is a structure map for the
left square.

Now let
$$
\Phi_{n}(\Xd)\co \map_{\CA/\bBL}(\Xd,\bEL{M}{n})\to
\map_{s\PAAlg/\tBL}(\piA \Xd,\tEL{M}{n})
$$
\noindent be the analogously defined map, with
\w[.]{\phi_{n}(X)=\pi_{0}\Phi_{n}(\Xd)}

Because \w{\piA} takes homotopy pushouts in \w{\CA} to homotopy
pushouts of simplicial \APa s, it follows that the source and target
of \w{\Phi_{n}(-)} take homotopy pushouts to homotopy pullbacks. Now
every object of \w{\CA} is, up to homotopy, a filtered colimit of
objects constructed from copies of \w{{A} \hotimes \bS{m}}  by finitely
many homotopy
pushouts. Thus, since source and target of \w{\Phi_{n}} take
filtered colimits to homotopy limits of simplicial sets, it suffices
to show that \w{\Phi_{n}({A} \hotimes \bS{m})} is a
$\pi_{0}$--equivalence
for all \w{m\geq 2} and \w[.]{A\in\calA} As \w{{A} \hotimes \bS{m}}
corepresents \w{\pin{n}(?)} in \w{\ho\CA/\bBL} and \w{\piA({A} \hotimes
\bS{m})}
corepresents \w{\pi_{n}\piA (?)} in \w{\ho s\PAAlg} for \w[,]{n\geq 2}
the result follows from the naturality of $\zeta$ and \fullref{deaem}.
\end{proof}

The restriction \w{n\geq 2} is needed because \w{\pi_{1}\piA (?)} is
not known to be corepresentable (see Dwyer, Kan and Stover
\cite[Section~7(ii)]{DKStoB}).

%
%
\begin{cor}\label{czero}
The functors \w{H^{\ast}_{\Lambda}(-;M)}  on
\w{\CA/\bBL} and \w{s\PAAlg/\tBL} are cohomology functors.
\end{cor}

\begin{proof}
This follows from Quillen \cite[Section~II.5]{QuiH}.
\end{proof}

\begin{remark}\label{rcoh}
If $\calC$ is the category \w[,]{\PAAlg} or more generally any
category of \Ta s as in \fullref{ttwo}, we have an equivalence
\[
H^{n}_{\Lambda}(\Gd;M) \cong
\pi_{0}\map_{s\RM{\Gd}/\BL}(\LL\Omega_{\Gd},\EL{M}{n}).
\]
Here \w{\LL\Omega_{\Gd}} denotes the \emph{cotangent complex}
associated to \w[,]{\Gd} defined by
\[
\LL\Omega_{\Gd} := \Omega_{\Gd'}\ast_{\Gd'}\Gd
\]
where \w{\Gd'} is a cofibrant replacement of \w{\Gd} in \w{\CA}
and the group of K\"{a}hler differentials \w{\Omega_{\Gd'}} is defined in
\eqref{Kaehler}.

\end{remark}

\begin{remark}\label{cohbyext}
In fact, this previous observation can	be carried a little
further. Given a (simplicial) \APa\ \w{\Gd} and a \w[-]{\Gd}module
$M$, define the {\it group of algebraic extensions}
\w{\exal{\Lambda}{\Gd}{M}} to be the set of equivalence classes of the
form \eqref{algext} with \w[.]{K = M} This set is a functor in both
variables (via pullbacks and pushouts) and forms an abelian group with
unit \w{M\ltimes \Gd} and addition induced by the diagonal
\w{\Gd\to\Gd\times_{\Lambda} \Gd} and the group addition
\w[.]{M\times_{\Lambda} M\to M}

Assume now that \w{\Gd} is cofibrant. Following Illusie \cite[III.1.2.3]{IlC},
there is a natural isomorphism
\begin{equation}\label{extiso1}
\exal{\Lambda}{\Gd}{\EL{M}{n}}\xra{\cong} H^{n+1}_{\Lambda}(\Gd;M)
\end{equation}

\noindent sending an algebraic extension \w{(\EL{M}{n} \to X \to \Gd)}
of simplicial \APa s to the induced homotopy coboundary
\w[.]{(\Gd \to\EL{M}{n+1})} For general \w[,]{\Gd} there is an isomorphism
\begin{equation}\label{extiso2}
H^{n+1}_{\Lambda}(\Gd;M)\cong\colim_{\mbox{Wk}(\Gd)}\,
\exal{\Lambda}{\Gd'}{\EL{M}{n+1}}
\end{equation}

\noindent where \w{\mbox{Wk}(\Gd)} is the category of cofibrant
replacements \w{\Gd' \to \Gd} in simplicial \APa s
\end{remark}

\subsection{The cohomology of a diagram}\label{scohdiag}

Let $\DD$ be a small category. Observe that a map of $\DD$--diagrams is
just a natural transformation: a collection of maps on objects which
commute with the maps in each diagram.
%

%
%
\begin{fact}\label{ftwo}
Given two functors \w[,]{X,Y\co \DD\to\calC} the set
\w{\Hom_{\calC^{\DD}}(X,Y)}
of diagram maps between them fits into the equalizer diagram
\begin{equation}\label{eqtwenty}
\Hom_{\calC^{\DD}}(X,Y)\hra\prod_{d\in\DD}\Hom_{\calC}(X_d,Y_d)\equaliz
\prod_{d,e\in\DD}\prod_{\eta\in\Hom_{\DD}(d,e)}\Hom_{\calC}(X_d,Y_e),
\end{equation}
\noindent where the two parallel arrows map to each factor indexed
by \w{\eta\co d\to e} in $\DD$ by the appropriate projection, followed by
\w[,]{Y(\eta)_{\ast}\co \Hom_{\calC}(X_d,Y_d)\to\Hom_{\calC}(X_d,Y_e)} or
\w[,]{X(\eta)^{\ast}\co \Hom_{\calC}(X_e,Y_e)\to\Hom_{\calC}(X_d,Y_e)}
respectively.
\end{fact}

\begin{remark}\label{rfib}
If $\calC$ is a simplicial model category, and \w{Y_{d}}
is an abelian group object for each \w[,]{d\in\Obj\DD} we can
replace the equalizer diagram \eqref{eqtwenty} by an exact sequence
of simplicial abelian mapping spaces (using the mapping
space construction of Quillen \cite[II.3.1]{QuiH})
%
\begin{equation}\label{eqtwentyone}
0\to\map_{\calC^{\DD}}(X,Y)\to\prod_{d\in\DD}
\map_{\calC}(X_d,Y_d)\xra{\xi}
\prod_{d,e\in\DD}\prod_{\eta\co d\to e}\map_{\calC}(X_d,Y_e),
\end{equation}
\noindent where $\xi$ is the difference of the two parallel arrows of
\eqref{eqtwenty}.
\end{remark}

If this were a fibration sequence after the mapping spaces are
restricted to appropriate over-categories, we could apply
\w{\pi_{0}} and compute cohomology in the diagram category directly from
the exact sequence.  However, it is not a fibration sequence in
general, so we concentrate for now on the special case of \w[.]{\DD=[1]}

\subsection{The cohomology of a map}\label{arrowcohom}

For the arrow category \w[,]{\MC} the exact sequence of
\eqref{eqtwentyone}, suitably
modified, is in fact a fibration sequence. To show this, we need some
technical results on model categories:
%

%
\begin{lemma}\label{fone}
Suppose
\mydiagram[\label{eqnineteen}]{
X \ar[dr]^{g} \rto^{f} & W  \dto^-{\psi} \\
	   & Z \\
}
\noindent is a diagram in a model category $\calC$
which commutes up to homotopy, with $X$ cofibrant and $\psi$ a
fibration.  There there is a homotopic map
\w{f \simeq f'\co X\to W} such that \w[.]{\psi\circ f'=g}
Dually, if
\mydiagram[\label{eqnineteenb}]{
X \ar[d]^{\phi} \ar[dr]^{f} \\
Y \ar[r]_{g}	      & Z \\
}
commutes up to homotopy, with $Z$ fibrant and $\phi$ a cofibration,
then there is a homotopic map \w{g \simeq g'\co Y\to Z} such that
\w[.]{f=g' \circ \phi}
\end{lemma}

\newcommand{\cyl}[1]{cyl({#1})}
\begin{proof}
Assume $\psi$ is a fibration. Cofibrancy of $X$ implies
\w{i_0\co X\to\cyl{X}} is an acyclic cofibration by Hirschhorn \cite[7.3.7]{HirM}.
Given a homotopy \w{H\co \cyl{X} \to Z} with
\w{H\circ i_{0}=\psi\circ f} and \w[,]{H\circ i_{1}=g} we
may use the left lifting property in
\begin{equation}{\label{ellp}}
\vcenter{\xymatrix@C=25pt{
X \ar[d]_{\text{acyc.~cof}}^{i_{0}} \ar[r]^{f} &
		 W \ar[d]^{\text{fib}}_{\psi} \\
{\cyl{X}} \ar@{.>}[ru]^{\hat{H}} \ar[r]^{H}	 & Z
}}
\end{equation}
\noindent to factor $H$ as \w[,]{\psi\circ\hat{H}} and set
\w[.]{f':=\hat{H}\circ i_{1}} If instead $\phi$ is a cofibration and
$Z$ is fibrant, use the dual argument.
\end{proof}

%
%
\begin{cor}\label{cone}
Suppose
\mydiagram[\label{eqnineteentwo}]{
X \dto^{\phi} \rto^{f} & W  \dto^-{\psi} \\
Y \rto^{g}	    & Z \\
}
\noindent is a commutative diagram in a model category $\calC$.
If $\psi$ is a fibration and $X$ is cofibrant, then to any homotopic
map \w{g' \simeq g} there corresponds a homotopic map \w{f' \simeq f}
such that \w[.]{\psi \circ f'=g'\circ \phi} Dually, if $\phi$ is a
cofibration and $Z$ is fibrant, then to any homotopic map \w{f'\simeq f}
there corresponds a homotopic map \w{g' \simeq g} such that
\w[.]{\psi\circ f'=g'\circ\phi}
\end{cor}

\begin{remark}\label{rhcm}
Since we assume that fibrations and weak equivalences in our
diagram categories are defined objectwise, then if $\phi$ is a
cofibrant  object in $\MC$ it follows that $\phi$ is a cofibration in
$\calC$ with cofibrant source. Thus if $\psi$ is a fibration
with fibrant target in $\calC$, it makes sense to consider homotopy
classes of maps
\w{[\phi,\psi]} in \eqref{eqnineteen} -- 
in fact, the mapping space \w{\map_{\MC}(\phi,\psi)} has homotopical
meaning, and \w[.]{[\phi,\psi]\cong\pi_{0}\map_{\MC}(\phi,\psi)}
\end{remark}

%
%
\begin{prop}\label{ptwo}
Let \w{\vartheta\co U\to V} be a fixed map in a simplicial model
  category $\calC$ and let \w{\phi\co X\to Y} and \w{\psi\co W\to Z} be maps
  in \w[.]{\MC/\vartheta} If $\phi$ is a cofibration with cofibrant
  source and \w{Z \to V} is a
fibration in $\calC$, with $W$ and $Z$ abelian group objects, then the
  restriction of the exact sequence of simplicial abelian mapping
  spaces from \fullref{rfib}
%
\begin{equation}\label{eqtwentytwo}
\map_{\MC/\vartheta}(\phi,\psi)\to
\map_{\calC/U}(X,W)\times\map_{\calC/V}(Y,Z)\xra{\xi}\map_{\calC/V}(X,Z)
\end{equation}
is a fibration sequence (in $\calS$).
\end{prop}

\begin{proof}
First, by Quillen \cite[Proposition~II.3.1]{QuiH}, we know that $\xi$ of
\eqref{eqtwentytwo}
is a fibration in $\calG$ (and so in $\calS$) if and only if it surjects
onto the basepoint component of the target space
\w{\map_{\calC/V}(X,Z)\in\calS} -- or equivalently, onto any component
of
$\map_{\calC/V}(X,Z)$ which it hits.

Now, if \w{k\co X\times\Delt{n}\to Z} is any map
in the image of $\xi$, then there are maps \w{f\co X\times\Delt{n}\to W}
in \w{\calC/U} and \w{g\co Y\times\Delt{n}\to Z} in \w{\calC/V} such that
in the (not commutative) diagram
%
\mydiagram[\label{eqtwentythree}]{
{X\otimes\Delt{n}} \ar[d]_{\phi \otimes \Id} \ar[rr]^{f} \ar[drr]^{k}
      && {W} \ar[d]^{\psi} \\
{Y\otimes\Delt{n}}  \ar[rr]_{g} && {Z}
}
\noindent we have \w{\psi\circ f-g\circ(\phi\otimes\Id)=k} in
\w[.]{\calC/V}

Finally, if \w{k'} is in the same component as $k$ in
\w[,]{\map_{\calC/V}(X,Z)}
we can write \w{\psi\circ f- g\circ(\phi\otimes\Id)\sim_{V} k'} (since
$X$ is cofibrant and $Z$ is fibrant in \w[)]{\calC/V} or
equivalently, since $\pm$ preserves homotopies, \w[,]{\psi\circ
   f-k'\sim_{V} g\circ(\phi\otimes\Id)}
where $\sim_{V}$ indicates homotopy in \w[.]{\calC/V}
By \fullref{fone} applied to the diagram
\mydiagram[\label{eqnineteenbagain}]{
{X\otimes\Delt{n}} \ar[d]^{\phi \otimes \Id} \ar[drr]^{\psi\circ f-k'} \\
{Y\otimes\Delt{n}} \ar[rr]_{g}		&& Z \\
}
viewed in \w[,]{\calC/V}
we can replace $g$ by a homotopic map $g'$
over $V$  such that
\w[.]{\psi\circ f-k'=g'\circ(\phi\otimes\Id)} But then \w[,]{\xi(f,g')=k'}
so $\xi$ indeed surjects onto the component of $k$.
\end{proof}

%
%
\begin{cor}\label{ctwo}
For \w[,]{\phi\co \Xd\to\Yd} a morphism in \w{\sC} over a map
\w[,]{B\lambda\co \BL_{0}\to\BL_{1}} suppose
\w{\psi\co E^{\Lambda_{0}}(M_{0},n)\to E^{\Lambda_{1}}(M_{1},n)}
is the morphism of extended Eilenberg--Mac\,Lane objects induced by a
module \w{\tau\co  M_{0}\to M_{1}} over
\w[.]{\lambda\co \Lambda_{0}\to\Lambda_{1}}
Then there is a long exact sequence
%
\begin{multline}\label{eqtwentyfour}
0\to H^{0}_{\lambda}(\phi,\tau) \to
H^{0}_{\Lambda_{0}}(\Xd;M_{0})\oplus H^{0}_{\Lambda_{1}}(\Yd;M_{1})
\xra{\psi_{\ast}-\phi^{\ast}} \\
H^{0}_{\Lambda_1}(\Xd;M_{1})\to
H^1(\phi,\tau)\to\cdots\to H^{n-1}_{\Lambda_{1}}(\Xd;M_{1})\to\\
H^{n}_{\lambda}(\phi;\tau)\xra{\theta}
H^{n}_{\Lambda_{0}}(\Xd;M_{0})\oplus H^{n}_{\Lambda_{1}}(\Yd;M_{1})
\xra{\psi_{\ast}-\phi^{\ast}}H^{n}_{\Lambda_{1}}(\Xd;M_{1})
\end{multline}
\noindent where $\theta$ is induced by the obvious forgetful functors.
\end{cor}

\begin{proof}
Recall from \fullref{rrss} that we may assume that our extended
Eilen\-berg--Mac\,Lane objects are strict abelian group objects, so that
the previous discussion applies. Note also that
\w{H^{n-r}_{\Gamma}(\Wd,N) \cong\pi_{r}\map_{s\calC}(\Wd,E^\Gamma(N,n))}
for \w[,]{\Wd\in\sC/B\Gamma} $N$ a $\Gamma$--module, and \w[.]{0\leq r
\leq n}
Similarly  \w[.]{H_{\lambda}^{n-r}(\phi,\tau)\cong
\pi_{r}\map_{\sC(\to)}(\phi,E^\lambda(\tau,n))}
Thus the fibration sequence \eqref{eqtwentytwo} yields the desired
long exact
sequence in homotopy (though the last map in \w{\pi_{0}} need not
be surjective).
\end{proof}

We can identify the image of \w{\psi_{\ast}-\phi^{\ast}} in
cohomological terms as
$$
\Ker(q_{\ast}\co H^{n}(\Xd;M_1)\to H^{n}(\Xd;C))\cap
\Image(\phi^{\ast}\co H^{n}(\Yd;M_1)\to H^{n}(\Xd;M_1)),
$$
where $q\co M_1\epic C:=\Cok(\tau)$\vsm.

\subsection{An example of the cohomology of a map}\label{scohmap}

Note that in the stable range any $\Lambda$--module is
\emph{trivial} -- that is, \w{\llra{-}{-}\equiv 0} (in the
notation of \fullref{rmod}) (although of
course it need not be trivial as an abelian \Pa -- that is,
compositions may be non-zero).

In our example, for \w{\Lambda:=\tru{n+2}\pis\bX} (\fullref{sspace}),
and \w[,]{M:=\Omega\Lambda} we have
\[
M_{i}=\begin{cases}
(\ZZ/2)\lra{\alpha} & \text{for }i=n-1\\
(\ZZ/2)\lra{\alpha\circ\eta} & \text{for }i=n\\
(\ZZ/4)\lra{\beta} & \text{for }i=n+1\\
0 & \text{for }i=n+2,
\end{cases}
\]
\noindent with $2\beta=\alpha\circ\eta^{2}$\vsm.

Since  \w{\PAn} is an abelian category, by the Dold--Kan
correspondence we can use chain-complex notation to describe a free
simplicial resolution \w{\eVd} of $\Lambda$ as follows:

%
%
{\small
\mydiagram[\label{fig1}]{
\gS{n+2}_{s} \ar[r]^{2} & \gS{n+2}_{t} \ar[rd]^{\eta} & & &
\gS{n+2}_{w} \ar[r]^{2} \ar[rdd]^{-\eta^{2}} &
\gS{n+2}_{y} \ar@{|->}[r] & \beta\\
   & & \gS{n+1}_{v} \ar[r]^{2} & \gS{n+1}_{u} \ar[rd]^{\eta} &
\amalg &\amalg & \\
   & & & & \gS{n}_{z} \ar[r]^{2} & \gS{n}_{x} \ar@{|->}[r] & \alpha\\
\eV_{5} \ar[r]^{\partial_{5}} & \eV_{4} \ar[r]^{\partial_{4}} &
\eV_{3} \ar[r]^{\partial_{3}} & \eV_{2} \ar[r]^{\partial_{2}} &
\eV_{1} \ar[r]^{\partial_{1}} & \eV_{0} \ar[r]^{\partial_{0}} &
\Lambda,
}}
\noindent (where \w[)]{\partial_{1}(w)=2y-x\circ\eta^{2}\in \eV_{0}} --
so we can calculate
$$
C^{\ast}:=\Hom_{\RM{\Lambda}}(\eVd,\Omega\Lambda)
$$
as follows
\begin{equation*}
\begin{array}{ccccccccccc}
C^{5} & \leftarrow 
& C^{4} & \leftarrow 
& C^{3} & \leftarrow 
& C^{2} & \leftarrow 
& C^{1} & \leftarrow 
& C^{0}\\
\| & & \| & & \| & & \| & & \| & & \|\\
0 & \xla{0} & 0 & \xla{0} & \ZZ/4 & \xla{2} & \ZZ/4 & \xla{2} & \ZZ/2 &
\xla{0} & \ZZ/2

\end{array}
\end{equation*}

\noindent which implies that
\[
H^{i}(\Lambda;\Omega\Lambda)=\begin{cases}
\ZZ/2& \text{for }i=0,3\\
0& \text{otherwise}\vsm.
\end{cases}
\]
Similarly, \w{\Hom(\eVd,\Omega\gS{n-1})} is
\w[,]{0\leftarrow 0\xla{0}\ZZ/24\xla{2}\ZZ/24\xla{12}\ZZ/2\xla{0}\ZZ/2}
so that
$$
H^{i}(\Lambda;\Omega\gS{n-1})=\begin{cases}
\ZZ/2& \text{for }i=0,3\\
0& \text{otherwise}
\end{cases}
$$
\noindent with
\w{\varphi_{\ast}\co H^{0}(\Lambda;\Omega\gS{n-1})\to
H^{0}(\Lambda;\Omega\Lambda)}
the identity, while
$$\varphi_{\ast}\co H^{3}(\Lambda;\Omega\gS{n-1})\to
H^{3}(\Lambda;\Omega\Lambda)
$$
is trivial (and similarly for $\psi$ of \fullref{rreal}).

On the other hand, since \w{\gS{n-1}} is a free \Pa, for any module
$M$ we have
\[
H^{i}(\gS{n-1};M)=\begin{cases}
M& \text{for }i=0\\
0& \text{otherwise.}
\end{cases}
\]
{From} the long exact sequence \eqref{eqtwentyfour} we conclude that:
%
\begin{equation}\label{eqthirtytwo}
H^{i}_{\varphi}(\varphi;\Omega\varphi)
=H^{i}_{\psi}(\psi;\Omega\psi)=
\begin{cases}
\ZZ/2& \text{for }i=3,4\\
0& \text{for } 0<i<3 \text{ or  }4< i.
\end{cases}
\end{equation}
%

%
%
\section{Realizations of a \APa}
\label{crpia}

Our aim now is to address the general realization question described
in the
introduction  --  namely,  given an $E^{2}$--model category
\w{\CA} and a \APa\ $\Lambda$, is there a {\it realization} of
$\Lambda$ in $\calC$ - that is, is there a $Y\in\calC$ such that
\w{\piA Y \cong \Lambda} as $\Pi$--algebras?

Before we state our main result, we need the following variation on the
Postnikov system:

\begin{defn}\label{dquap}
A \emph{quasi-Postnikov tower} for an \APa\ $\Lambda$ is a tower of
fibrations
%
\begin{equation}\label{eseven}
\dotsb\xra{p\q{n+1}}\Xn{n+1}\xra{p\q{n}}\Xn{n}\xra{p\q{n-1}}\dotsb\xra{p\q{0}}
\Xn{0}\simeq\bBL
\end{equation}
\noindent in \w{\sC} such that \w{\piA\Xn{n}\simeq\tEL{\OA{n+1}}{n+2}}
for every \w[,]{n>0} with the sections \w{s\co \tBL\to\piA\Xn{n}}
(\fullref{rsec})
induced by the maps \w[.]{p\q{n}}
The object \w{\Xn{n}\in\sC} will be called an
\emph{$n$th quasi-Postnikov section} for $\Lambda$.
\end{defn}

\begin{remark}\label{rqpost}
Thus a tower \eqref{eseven} is a quasi-Postnikov tower for $\Lambda$ if
%
\begin{equation}\label{eqthirteen}
\pi_{k}\piA\Xn{n}\cong\begin{cases}
		\Lambda  &\text{for } k=0,\\
		\OA{n+1} &\text{for } k=n+2,\\
		 0	 &\text{otherwise},\end{cases}
\end{equation}
\noindent and it is equipped with maps \w{\ruq{n}\co \tBL\to\piA\Xn{n}} over
\w[,]{\tBL} for each \w[,]{n\geq 0} commuting with the maps
\w[.]{p\q{n}_{\#}}

We then deduce from the exact sequence \eqref{eqfour} that
%
\begin{equation}\label{eqtwelve}
\pin{k}{\Xn{n}}\cong\begin{cases}
		\OA{k} &\text{for } 0\leq k\leq n,\\
		 0     &\text{otherwise}.\end{cases}
\end{equation}
Note that \eqref{eqtwelve} implies in turn that the (ordinary)
Postnikov
sections \w{\bP{k}\Xn{n}} of \w{\Xn{n}} constitute quasi-Postnikov
sections for $\Lambda$, for \w{k\leq n} (see also Blanc, Dwyer and Goerss
\cite[Proposition~9.11]{BDGoerR}).
\end{remark}

We are now in a position to state the two key results addressing our
realization question (the proofs are deferred to
Sections~\ref{srp}--\ref{srpc}):

\begin{thm}\label{tfour}
If \w{\CA} is an $E^{2}$--model category and \w[,]{\Lambda \in \PAAlg}
the following are equivalent:
\begin{enumerate}
\renewcommand{\labelenumi}{(\arabic{enumi})}
\item $\Lambda$ is realizable -- that is, there is a \w{Y\in\calC}
with
  \w[;]{\piA Y\cong\Lambda}
\item There is an \w{\Xd\in\sC} with \w[.]{\piA\Xd\simeq \tBL}
\item There is a quasi-Postnikov tower for $\Lambda$.
\end{enumerate}
\end{thm}

\begin{thm}\label{tclass}
Let \w{\Xn{n-1}\in\sC} be an $(n{-}1)$st quasi-Postnikov section for
a \APa\ $\Lambda$. Then:
\begin{enumerate}
\renewcommand{\labelenumi}{(\alph{enumi})}
\item Up to homotopy, there is a unique \w{\Xn{n}\in\sC} satisfying
  \eqref{eqthirteen} and \eqref{eqtwelve}, with
  \w[.]{\bP{n-1}\Xn{n}=\Xn{n-1}}
\item  This \w{\Xn{n}} is an $n$th quasi-Postnikov section for
  $\Lambda$ if and only if the $(n{+}2)$nd $\tk{}$--invariant for
  \w{\piA\Xn{n}} vanishes in \w[.]{\HL{n+3}{\tBL}{\OA{n+1}}}
\item In that case, \w{\Xn{n+1}} exists, by (a); furthermore, the
  different choices for the map  \w{p\q{n}\co \Xn{n+1}\to\Xn{n}} -- or
  equivalently, choices of the	section
  \w{\ts{n}\co \tBL\to\tEL{\OA{n+1}}{n+2}=\piA\Xn{n}} of \fullref{rsec} -- are in one-to-one correspondence with elements
  of \w[.]{\HL{n+2}{\tBL}{\OA{n+1}}}
\end{enumerate}
\end{thm}
\noindent Compare Baues \cite[Chapter~D, (7.9)]{BauCF}.

Our approach to constructing an \w{\Xd} in \fullref{tfour}~(2)
will be inductive, using its Postnikov system, which serves as a
quasi-Postnikov tower for $\Lambda$. Thus at each stage we will have
the obstruction of \fullref{tclass}~(b) to moving up one more level.
To explain why this works (and prove the two Theorems), we shall need
some facts about:

\subsection{Connections between the Postnikov systems}\label{ccps}

Given any simplicial object \w[,]{\Xd\in\sC} consider its $n$th Postnikov
section \w[,]{\bP{n}\Xd} for some \w[,]{n>0} and let
\w[.]{\As\DEF\pin{0}{\Xd}=\pi_{0}\piA\Xd}
We want to describe the simplicial \APa\ \w{\piA\bP{n}\Xd} (up to
homotopy) in terms of \w[,]{\piA\Xd} and whatever other information
is necessary.

First, observe that \eqref{eqfour} also implies
%
\begin{equation}\label{efive}
\pi_{k}\piA\bP{n}\Xd\cong\begin{cases}
	 \pi_{k}\piA\Xd & \text{for } k\leq n,\\
	 \Cok(h^{X}_{n+1}\co \pin{n+1}{\Xd}\to\pi_{n+1}\piA{\Xd}) & \text{for }
	   k=n+1,\\
	 \Omega\pin{n}{\Xd} & \text{for }k=n+2.\\
0 & \text{otherwise}.
\end{cases}
\end{equation}
In particular, when \w[,]{\piA\Xd\simeq\tBL} \eqref{efive} simplifies to
%
\begin{equation}\label{enine}
\pi_{k}\piA\bP{n}\Xd\cong\begin{cases}
	 \As & \text{for } k=0,\\
	 \Omega^{n+1}\As & \text{for }k=n+2,\\
0 & \text{otherwise},
\end{cases}
\end{equation}
%
%
\begin{lemma}\label{lthree}
For any \w[,]{\Xd\in\sC} we have a homotopy fibration sequence
$$
\piA\bP{n+1}\Xd\xra{p\q{n}_{\#}}\piA\bP{n}\Xd\xra{(\bk{n})_{\#}}
\piA\bE{\As}{\pin{n+1}{\Xd}}{n+2}
$$
in \w{s\PAAlg/\tBL}
(that is, a homotopy pullback square over \w[).]{\tBL}
\end{lemma}

\begin{proof}
\fullref{srmc}(b) implies that
$$
(\bk{n})_{\#}\co \piA\bP{n}\Xd\to\piA\bE{\As}{\pin{n+1}{\Xd}}{n+2}
$$
is an $\calA$--fibration over \w[.]{\piA\bBL} Denote its fiber by
\w[,]{\Fd}
with a natural map of simplicial \APa s
\w[.]{\varphi\co \piA\bP{n+1}\Xd\to\Fd}

Because the functors \w{\pi_{k}\piA\co \sC\to\PAAlg} are corepresentable for
\w{k>1} (cf Dwyer, Kan and Stover \cite[Section~7.4]{DKStoB}), applying
\w{\piA} to the
homotopy pull-back \eqref{eqfive} yields a ``quasi-fibration'' of
simplicial \APa s, and so a long exact sequence in homotopy (in
dimensions $\geq 2$), which implies that \w{\varphi_{\#}} is an
isomorphism in dimensions $\geq 2$; since this is trivially true in
dimensions $0$ and $1$, $\varphi$ is a weak equivalence.
\end{proof}

%
%
\begin{lemma}\label{fthree}
If we write \w[,]{\Ed:=\bEL{\pin{n+1}{\Xd}}{n+2}} then applying
\w{\pi_{n+2}\piA} to the $k$--invariant \w{\bk{n}\co \bP{n}\Xd\to\Ed}
yields the homomorphism \w{s_{n+1}\co \Omega\pin{n}{\Xd}\to\pin{n+1}{\Xd}}
of \eqref{eqfour}.
\end{lemma}

\begin{proof}
First, note that, in the commutative diagram
\begin{small}
$$
\xymatrix
{
\quad & 0 \ar[r] \ar[d] & \pin{n+1}{\Omega\Ed} \ar[d]^{\cong}
\ar[r]^{\cong} & \pi_{n+1}\piA\Omega\Ed \ar[d] \\
  \pi_{n+2}\piA\bP{n+1}\Xd \ar[d]\ar[r]^{\partial_{n+2}} &
\Omega\pin{n}{\bP{n+1}\Xd} \ar[d]^{\cong} \ar[r]^{s_{n+1}} &
\pin{n+1}{\bP{n+1}\Xd} \ar[d] \ar[r]^{h_{n+1}} &
\pi_{n+1}\piA\bP{n+1}\Xd \\
  \pi_{n+2}\piA\bP{n}\Xd \ar[d]_{(\bk{n})_{\#}} \ar[r]^{\cong} &
\Omega\pin{n}{\bP{n}\Xd} \ar[r] & \pin{n+1}{\bP{n}\Xd}=0 & \\
  \pi_{n+2}\piA{\Ed} & &
}
$$
\end{small}%
the isomorphisms of \w{\pi_{n+2}\piA\bP{n}\Xd} with
\w[,]{\Omega\pin{n}{\bP{n+1}\Xd}} and \w{\pin{n+1}{\bP{n+1}\Xd}} with
\w[,]{\pi_{n+1}\piA\Omega\Ed} are natural.  Also, the columns here are
exact either by the long exact sequence in \w{\pinat_{\ast}} for a
fibration in \w[,]{\sC} or by \fullref{lthree}.

The result now follows from the naturality of the exact sequence
\eqref{eqfour},
applied to the fibration sequence
$$\Omega\Ed\simeq\bEL{\pin{n+1}{\Xd}}{n+1}\to\bP{n+1}\Xd\to\bP{n}\Xd
\xra{\bk{n}} \Ed.\proved$$
\end{proof}

%
%

\begin{lemma}\label{scalc}
If \w[,]{\piA\Xd\simeq\tBL} then the spiral exact sequence
\eqref{eqfour} for \w{\Xd} from \w{\pi_{n+3}\piA\Xd} down is
determined by the  homomorphism
$$\w[.]{\partial_{n+3}^{\star}\co\pi_{n+3}\piA\Xd\to\Omega\pin{n+1}{\Xd}}$$
\end{lemma}

\begin{proof}
First,
observe that given \w[,]{\bP{n}\Xd} we know the exact sequence
\eqref{eqfour}
for \w{\Xd} only from \w{\Omega\pin{n-1}{\Xd}} down. However, when
\w{r\q{n}_{\#}\co \pis\piA\Xd\to\pis\piA\bP{n}\Xd} is also known, and
\w[,]{\piA\Xd\simeq\tBL} then all we need in order to determine
\eqref{eqfour} for \w{\Xd} from \w{\pi_{n+3}\piA\Xd} down is the
homomorphism
\w{(r\q{n+1}_{\#})_{n+3}\co \pi_{n+3}\piA\Xd\to\pi_{n+3}\piA\bP{n+1}\Xd}
-- which
is just
\w[.]{\partial_{n+3}^{\star}\co \pi_{n+3}\piA\Xd\to\Omega\pin{n+1}{\Xd}}
\end{proof}

\begin{lemma}\label{ffive}
If \w{\tk{n+1}(\piA\Xd)\co \tP{n+1}\piA\Xd\to\tEL{\pi_{n+2}\piA\Xd}{n+3}} is
the $(n{+}1)$st  $\tk{}$--invariant for \w[,]{\piA\Xd} then the
$(n{+}1)$st $\tk{}$--invariant
$$
\tk{n+1}(\piA\bP{n}\Xd)\co \tP{n+1}\piA\bP{n}\Xd\to\tEL{\Omega\pin{n}{\Xd}}{n+3}
$$
\noindent satisfies
\w[.]{(\partial^{\star}_{n+2})_{\ast}\circ \tk{n+1}(\piA\Xd) =
\tk{n+1}(\piA\bP{n}\Xd) \circ \tP{n+1}(r\q{n}_{\#})}
\end{lemma}

\begin{proof}
This follows from the naturality of the $\tk{}$--invariants (Ax~\ref{ax3} of
\fullref{demc}) and \fullref{scalc}.
\end{proof}

\subsection[Proof of Theorem~\ref{tfour}]{Proof of \fullref{tfour}}\label{srp}

\noindent (1) $\Longleftrightarrow$ (2)\qua Given $Y$, let
\w[.]{\Xd:=\cons{Y}}
Conversely, if \w{\Xd\in \sC/\bBL} satisfies \w[,]{\piA\Xd\simeq
\tBL} then
by Ax~\ref{ax4} of \fullref{demc}, there is a functor \w{J\co  \CA \to \calC}
equipped with an isomorphism
$$
[A,J\Xd]_{\calC}\cong\Hom_{\PAAlg}(\piA A,\Lambda),
$$
\noindent natural in \w[.]{A\in \calA} Thus \w{\piA J\Xd \cong
  \Lambda} as \APa s, by Yoneda's Lemma, so we can take \w[\vsm.]{Y:=J\Xd}

\noindent (2) $\Longleftrightarrow$ (3)\qua By Blanc, Dwyer and Goerss
\cite[Proposition~9.11]{BDGoerR}
we know that \w{\piA\Xd\simeq\tBL} if and only if
\w[.]{\piA\bP{n}\Xd\simeq \tEL{\OA{n+1}}{n+2}}

Thus given \w{\Xd} with \w[,]{\piA\Xd\simeq \tBL} the ordinary Postnikov
tower \w{\bP{k}\Xd} of \w{\Xn{n}} constitutes a quasi-Postnikov
tower for $\Lambda$, by \eqref{enine}.

Conversely, given a quasi-Postnikov tower \eqref{eseven} for
$\Lambda$, let \w[.]{\Xd \DEF \holim_n \Xn{n}}	Since
\w{\tP{n+1}\ruq{n}\co \tBL\to\tP{n+1}\piA \Xn{n}} is a weak equivalence
for each $n$, the maps \w{\ruq{n}} induce a weak equivalence
\w[.]{r\co \tBL\xra{\simeq}\piA\Xd}
\hfill$\Box$

\subsection[Proof of Theorem~\ref{tclass}]{Proof of \fullref{tclass}}\label{srpc}

Let \w{\Xn{n-1}} be an $(n{-}1)$st quasi-Postnikov section for $\Lambda$.
By assumption
$$\w[,]{\piA\Xn{n-1}\simeq\tEL{\OA{n}}{n+1}}$$
and the map
\w{\ruq{n-1}\co \tBL\to\piA\Xn{n-1}} is the required section.
\begin{enumerate}
\renewcommand{\labelenumi}{(\alph{enumi})}
\item In order to construct \w[,]{\Xn{n}} we must choose a suitable
$(n{-}1)$st $\bk{}$--invariant
\w[.]{\bk{n-1}\in[\Xn{n-1},\bE{\As}{\OA{n}}{n+1}]_{\bBL}}
Note that using the long exact sequence in \w{\pinat} for a
fibration over \w[,]{\bBL} combined with \eqref{eqfour},
automatically ensures that any such choice yields \w{\Xn{n}}
satisfying \eqref{eqthirteen} and \eqref{eqtwelve}.

We can use the map
\w{\zeta\co \piA\bEL{\OA{n}}{n+1}\to\tEL{\OA{n}}{n+1}} of \fullref{psix}
to define \w{\bk{n-1}\co \Xn{n-1}\to\bEL{\OA{n}}{n+1}}
(uniquely up to homotopy) by specifying
$$
\zeta\circ(\bk{n-1})_{\#}\co \piA\Xn{n-1}\to\tEL{\OA{n}}{n+1}.
$$
\noindent Since \w[,]{\piA\Xn{n-1}\simeq\tEL{\OA{n}}{n+1}} the
functoriality of Ax~\ref{ax2} of \fullref{demc} implies that such a map is
uniquely determined up to homotopy by a map of $\Lambda$--modules
\w[,]{\varphi\co \OA{n}\to\OA{n}} and by \fullref{fthree} this
$\varphi$ must be the given isomorphism
\w[,]{(s_{n+1})\co \Omega\pin{n-1}{\Xn{n-1}}\to\OA{n}} if the
quasi-Postnikov tower we are constructing for $\Lambda$ is to be a
Postnikov tower in \w[.]{\sC} (Note that by \fullref{scalc}, we
already know the long exact sequence \eqref{eqfour} for \w{\Xn{n}}
from \w{s_{n+1}} down.) Thus the candidate for \w{\Xn{n}} over
\w[,]{\Xn{n-1}}
satisfying \eqref{eqthirteen} and \eqref{eqtwelve}, is determined
uniquely
up to homotopy by \w[.]{\Xn{n-1}}
\item There is only one possible obstruction to \w{\Xn{n}} (the homotopy
fiber
  of \w{\bk{n-1}} in \w[)]{\sC/\bBL} being an $n$th quasi-Postnikov
  section for $\Lambda$: the non-existence of the lift
\w[.]{\ruq{n}\co \tBL\to\piA\Xn{n}} However, since
\w[,]{\tP{n+1}\piA\Xn{n}\simeq\tBL} by \eqref{eqtwelve},
we may use the long exact sequence in \w{\piA} for the fibration sequence
\begin{equation}\label{eqfifteen}
\piA\Xn{n}=\tP{n{+}2}\piA\Xn{n}{\xra{\tilde{p}\q{n{+}2}}}\tP{n{+}1}\piA\Xn{n}
{\xra{\tk{n{+}1}}}\tEL{\OA{n{+}1}}{n{+}3}
\end{equation}
\noindent over \w{\tBL} to deduce that \w{\ruq{n-1}} lifts to \w{\ruq{n}}
if and only if \w{\tk{n+1}} is null in \w[.]{s\PAAlg/\tBL}

More precisely, we want $\ruq{n}$ to map to the homotopy pullback
(Ax~\ref{ax3} of \fullref{demc}) in
%
\begin{equation}
\label{eqfourteen}
\xymatrix{{\tP{n+1}\tBL} \ar@/^1pc/[drr]^{\simeq} \ar@/_1pc/[ddr]_{\simeq}
\ar@{.>}[dr]|-{\ruq{n}} \\ & \piA\bP{n}\Xd \rto \dto & \tBL
\dto^-{\tk{n+1}} \\
& \tBL	 \rto^-{\tilde{s}}  & \tEL{\OA{n+1}}{n+3},}
\end{equation}
\noindent which is possible if and only if \w{\tk{n+1}} is homotopic
to the given homotopy section \w[.]{\tilde{s}\co \tBL\to
\tEL{\OA{n+1}}{n+3}}
\item Since the fiber (over \w[)]{\tBL} of \w{\tilde{p}\q{n+2}}
in \eqref{eqfifteen} is \w[,]{\tEL{\OA{n+1}}{n+2}} the possible choices
for such lifts are distinguished by elements of
$$
[\tBL,\tEL{\OA{n+1}}{n+2}]_{\tBL}=H^{n+2}(\tBL/\As,\OA{n+1}),
$$
which are just choices for
$\w{\partial_{n+3}^{\star}\co\pi_{n+3}\tBL\to\Omega\pin{n+1}{\Xn{n{+}1}}}$
(see \fullref{scalc}). These determine the identification of
\w{\piA\Xn{n}} with \w[,]{\tEL{\OA{n+1}}{n{+}2}} which is the only freedom
in the
inductive procedure we have described.
\end{enumerate}

\begin{remark}\label{stage1}

To appreciate the explicit inductive construction of these
obstructions provided in the above proof, let us examine more carefully
the
first step in realizing a \APa\ $\Lambda$:

Note first that, from the spiral exact sequence and Postnikov
sections, the homotopy groups of \w{\bBL} fit into the algebraic
extension
$$
\tEL{\Omega \Lambda}{2} \to \pis\bBL \to \tBL,
$$
\noindent and so yields an element of
\w{\exal{\Lambda}{\tBL}{\tEL{\Omega\Lambda}{2}}} (see \fullref{cohbyext}).
Using \eqref{extiso1}, we may view this extension as an element
of \w[,]{H^{3}(\tBL/\As,\Omega\Lambda)} which is precisely the
first obstruction to realizing $\Lambda$. Note that by Ax~\ref{ax4} of
\fullref{demc}, this obstruction is natural in $\Lambda$. See Benson, Krause
and Schwede \cite{BKS}
for a similar perspective on the obstructions to realizing modules
over the Tate cohomology of a group $G$ as the group cohomology of a
$G$--module.
\end{remark}

\begin{remark}\label{rapplic}
The realization problem, as formulated in this section, and its
solution in \fullref{tfour} applies to \Pa s  associated to any of
the categories listed in \fullref{stcat} - $n$--connected spaces,
$p$--local or
rational spaces, $n$--types (and so on) -  as well as any diagrams of
such \Pa s. Note, however, that realization is a tautology when
$\calC$ itself had a trivial model category structure -- e.g., if
$\calC=\Alg{\Theta}$ is a variety of universal algebras.
\end{remark}

%
%
\section{Realizing maps of \Pa s}
\label{crfpi}

We now examine the diagram realization question in more detail for the
simplest non-trivial case: a single map of (ordinary) \Pa s
\w[.]{\varphi\co \Lambda\to\Gamma}

\subsection{Maps of realizable \Pa s}\label{smrp}

Assume for simplicity that the \Pa s $\Lambda$ and $\Gamma$ are
realizable, and replace them by cofibrant simplicial models
\w{\psi\co \Kd\to\Ld} in \w[.]{s\PAlg}

Note that if we are \emph{given} realizations \w[,]{\Vd} \w{\Wd} for
\w{\Kd} and \w[,]{\Ld} respectively (equivalently: for $\Lambda$ and
$\Gamma$), we have the usual obstruction theory for lifting
\w{f^{0}:=\bB\phi\circ\bp{0}\co\Vd\to\bB\Gamma=\bP{0}\Wd}
through the successive Postnikov stages for \w[,]{\Wd} with
the existence and difference obstructions all lying in the Quillen
cohomology groups \w[.]{H^{\ast}(\Vd/\bB\Gamma;\Omega^{n}\Gamma)}
However, in our approach we want to choose the realizations for the
\Pa s\ $\Lambda$ and $\Gamma$, and for the map $\varphi$,
simultaneously -- again by induction on the quasi-Postnikov system.

At the $n$th stage, we assume that we have a map of simplicial spaces
\w[,]{\fn{n}\co \Xn{n}\to\Yn{n}} where:

\begin{enumerate}
\renewcommand{\labelenumi}{\alph{enumi})}
\item $\Xn{n}\simeq\bP{n}\Xn{n}$ and $\Yn{n}\simeq\bP{n}\Yn{n}$; and
\item $\tP{n}(\fn{n})_{\#}\co \tP{n}\piA\Xn{n}\to\tP{n}\piA\Yn{n}$ is
$\varphi_{\ast}\co \tBL\to\tB\Gamma$.
\end{enumerate}

Our goal is to extend $f$ to $(n{+}1)$--stage Postnikov pieces. Because
the sections $\ts{n}^{\Lambda}\co \tBL\to\tEL{\OA{n+1}}{n+2}$ and
$\ts{n}^{\Gamma}\co \tB\Gamma\to\tE{\Gamma}{\Omega^{n+1}\Gamma}{n+2}$
will ultimately be induced by the natural Postnikov maps
\w{\Wd\to\bP{n}\Wd\simeq\Xn{n}} and \w[,]{\Vd\to\bP{n}\Vd\simeq\Yn{n}}
say, we know that if \w{\fn{n}} extends we will have naturality for the
sections, so our object is to choose \w{\ts{n}^{\Lambda}} and
\w{\ts{n}^{\Gamma}} so that the diagram
%
\mydiagram[\label{eqeighteen}]{
\tBL \rto^{\varphi_{\#}} \dto^{\ts{n}^{\Lambda}} &
\tB\Gamma\dto^-{\ts{n}^{\Gamma}}\\
\tEL{\OA{n+1}}{n+2}\simeq\piA\Xd \rto^-{f_{\#}} &
\tE{\Gamma}{\Omega^{n+1}\Gamma}{n+2} \\
}
\noindent commutes up to homotopy. This means that
$(\ts{n}^{\Lambda},\ts{n}^{\Gamma})$ is just the obstruction class
in $H^{n+2}_{\varphi}(\varphi;\Omega^{n}\varphi)$ described by
\fullref{tfour}.

\subsection{An example of the obstructions to realizability}\label{sobst}

We now apply the above theory to the map of \Pa s
$\varphi\co \Lambda\to\gS{n-1}$ considered in \fullref{scohmap}.
By \cite[Theorem~3.16]{BlaCW}, we know that the resolution
\eqref{fig1}, as well as the constant free resolution
\w[,]{\eWd\to\gS{n-1}}
are realizable by simplicial spaces.

The relevant part of the realization of \eqref{fig1} is
described in \eqref{fig2}, where the indexing is based on the
Stover resolution comonad in the obvious way,
with \w{d_{0}} on \w{\bS{n+2}_{\lra{\beta,2}-\lra{\alpha,\eta^{2}}}}
equal to the difference of the degree $2$ map to \w{\bS{n+2}_{\beta}} and
\w{\eta^{2}} to \w[,]{\bS{n}_{\alpha}}, and all face maps \w{d_{1}} and
\w{d_{2}} are inclusions.

%
%
\begin{figure}[ht!]
{\small
\begin{equation}{\label{fig2}}
\vcenter{\xymatrix@C=11pt{
& \bS{n+2}_{\lra{\beta,2}-\lra{\alpha,\eta^{2}}}
\ar[rrdd]^>>>>>>>>>>>>>>>>>>>>>>>>>{d_{0}=2} \ar[rrd]^{d_{1}} & & & \\
\bS{n+1}_{\lra{\alpha,2,\eta}} \ar[r]^-{d_{2}} \ar[rd]^{d_{1}}
\ar[rdd]_{d_{0}=\eta} &
\mbox{$\bS{n+1}_{\lra{\alpha 2,\eta}}\cup\be{n+2}_{G,C\eta}$}
\ar[rrd]^>>>>>>>>>>>>>>>>>>>>>>>>{d_{1}}
\ar@{.>}[rrdd]^>>>>>>{d_{0}=C\eta} &
\ar@{.>}[rddd]^>>>>>>>>>>>>>>>>>>>>>>>>>>>>>>>>>>>>>>>>>
{\mbox{\hspace*{-2.5mm}{\scriptsize $-\eta^{2}$}}} &
\mbox{\hspace*{4mm}
$\bS{n+2}_{\lra{\beta 2}-\lra{\alpha\eta^{2}}}\cup\be{n+3}_{H}$}
\ar@{.>}[rddd]^>>>>>>>>>>>>>>>>>>>>{\varepsilon=H} & \\
   & \mbox{\hspace*{10mm}$\bS{n+1}_{\lra{\alpha,2\eta}}
\cup\be{n+2}_{\alpha,F}$}
\ar@{->}[rr]_<<<<<<<{d_{1}} \ar@{.>}[rrdd]^>>>>>>>>>>>>>{d_{0}=F} & &
\mbox{\raisebox{-2.2ex}{\hspace*{2mm}$\bS{n+2}_{\beta}=\bS{n+1}_{\alpha2\eta}
\cup\be{n+2}_{\alpha\circ F}\cup\be{n+2}_{G\circ C\eta}$}}
\ar[rdd]_<<<<{\mbox{\hspace*{-2mm}{\scriptsize $\varepsilon=\beta$}}} & \\
   & \bS{n}_{\lra{\alpha,2}} \ar[rr]^<<<<<<<{d_{1}} \ar[rrd]_{d_{0}=2} & &
{\mbox{\hspace*{3mm}$\bS{n}_{\alpha 2}\cup\be{n+1}_{G}$} }
\ar[rd]_<<<<<<{\varepsilon=G} & \\
  &  & & \bS{n}_{\alpha} \ar[r]^{\varepsilon=\alpha} & \bX\\
\bV_{2} \ar[r] \ar@<2.5mm>[r] \ar@<-2.5mm>[r]  &
\bV_{1} \ar@<1mm>[rr] \ar@<-1mm>[rr] &	& \bV_{0} \ar[r] & \bX,
}}
\end{equation}
\caption{A minimal free resolution $\bVd$ in $s\Top$}
}
\end{figure}

The inductive approach to realizing \w{\varphi\co \Lambda\to\gS{n-1}}
described in \fullref{smrp} begins with \w[,]{\fn{0}\co \Xn{0}\to\Yn{0}}
which is just \w[.]{\bB\varphi\co \bBL\to\bB\gS{n-1}} Moreover, the
proof of \fullref{tfour} shows that this always extends uniquely
to \w{\fn{1}\co \Xn{1}\to\Yn{1}} (although the lifting \w{\ruq{1}}
as required in \fullref{rqpost} need not exist).

The construction of Postnikov systems (Ax~\ref{ax1} of Theorems~\ref{ttwo}
and~\ref{tfive}) shows that the existence of \w{\fn{1}}
is equivalent to having a $2$--truncated augmented simplicial space
\w{\bVd'\to\bS{n-1}} realizing the augmented simplicial \Pa\
\w{\eVd\to\gS{n-1}}
induced by \w[.]{\varphi\co \Lambda\to\gS{n-1}}

Using \fullref{fone}, we may assume that the composite of the maps
$$
\bS{n}\xra{2}\bS{n}\xra{\eta}\bS{n-1}
$$
\noindent is actually null, so we can describe \w{\bVd'} explicitly by
\eqref{fig3}. Moreover, \w[,]{\Xn{1}} and thus  \w[,]{\bVd'} is
unique up to homotopy (in \w[).]{s\Top}

%
%
\begin{figure}[ht!]
{\small
\begin{equation}{\label{fig3}}
\vcenter{\xymatrix@C=14pt{
& \bS{n+2}_{\lra{2\nu,2}-\lra{\eta,\eta^{2}}}
\ar[rrdd]^>>>>>>>>>>>>>>>>>>>>>>>{d_{0}=2} \ar[rrd]^{d_{1}=\incl.}  &
& &\\
\bS{n+1}_{\lra{\eta,2,\eta}} \ar[r]^-{d_{2}=0} \ar[rd]^{d_{1}=\incl.}
\ar[rdd]_{d_{0}=\eta} & \ast &
\ar@{.>}[rddd]^>>>>>>>>>>>>>>>>>>>>>>>>>>>>>>>>>
{\mbox{\hspace*{-2mm}{\scriptsize $-\eta^{2}$}}} &
\mbox{\hspace*{4mm}
$\bS{n+2}_{\lra{4\nu}-\lra{\eta^{3}}}\cup\be{n+3}_{H}$}
\ar@{->}[rddd]^>>>>>>>>>>>>>>>>>>>>{\varepsilon=H} & \\
   & \mbox{\hspace*{10mm}$\bS{n+1}_{\lra{\eta,2\eta}}
   \cup\be{n+2}_{\eta,F}$}
\ar@{->}[rr]^<<<<<<<{d_{1}=0}
\ar@{.>}[rrdd]^>>>>>>>>>>>>>>>>>>>>>>>>>>>>>>{d_{0}=F} & &
\bS{n+2}_{2\nu} \ar[rdd]_<<<<{\varepsilon=2\nu} & \\
  & \bS{n}_{\lra{\eta,2}} \ar[rr]^<<<<<<<{d_{1}=0} \ar[rrd]_{d_{0}=2}
& &\ast& \\
  &  & & \bS{n}_{\eta} \ar[r]^{\varepsilon=\eta} & \bS{n-1}\\
\bV'_{2} \ar[r] \ar@<2.5mm>[r] \ar@<-2.5mm>[r]	&
\bV'_{1} \ar@<1mm>[rr] \ar@<-1mm>[rr] &  & \bV'_{0} \ar[r] & \bS{n-1},
}}
\end{equation}
\caption{An augmentation of \w{\bVd'} to $\bS{n-1}$}
}
\end{figure}

However, in constructing \w{\bVd'\to\bS{n-1}} we have ``distorted'' the
original augmented simplicial space \w{\bVd\to\bX} in such a way
that we no longer have a strict augmentation \w[.]{\bVd'\to\bX}

We can see this geometrically, using the Toda bracket
%
\begin{equation}\label{eqthirty}
\lra{\eta,2,\alpha}=\{\beta,\beta+\alpha\circ\eta^{2}\}
\subseteq\pi_{n+2}\bX
\end{equation}
\noindent (see, for example, \cite[Section~6]{BlaH}), which we used in the
decomposition
$$
\bS{n+2}_{\beta}=
\bS{n+1}_{\alpha2\eta}\cup\be{n+2}_{\alpha\circ F}\cup\be{n+2}_{G\circ
  C\eta}
$$
in \eqref{fig2}. Because we no longer have this in \eqref{fig3},
we must have \w{0\in\lra{\eta,2,\alpha}} for any augmentation
\w{\alpha\co \bS{n}\to\bX} on \w{\bS{n}\subseteq\bV'_{0}}

More formally, \eqref{eqthirty} yields a non-vanishing second-order
homotopy operation in \w{[\Sigma\bV'_{2},\bX]} which is the
obstruction to rectifying the homotopy augmentation \w{\bVd'\to\bX}
realizing \w[,]{\eVd\to\Lambda} using
\cite[Theorem~7.13 and Lemma~5.12]{BlaHH}. But then we may use the
equivalent obstruction theory of Blanc, Dwyer and Goerss
\cite{BlaAI,BDGoerR} to deduce that
the $\tk{}$--invariant \w{\tk{1}\in\HL{3}{\tBL}{\OA{}}\cong\ZZ/2}
does not
vanish, for the choice of \w{\Xn{0}} described in \eqref{fig3}
(with \w{\eta\co \bS{n}\to\bS{n-1}} replaced by \w{\alpha\co \bS{n}\to\bX}
and \w{2\nu} replaced by \w[).]{\beta\co \bS{n+2}\to\bX}

However, since the $\tk{}$--invariants are natural (\fullref{dki}), we deduce from the long exact sequence
\eqref{eqtwentyfour}
that the corresponding obstruction for the diagram -- that is,
\w{\tk{1}\in H^{3}(\varphi;\Omega\varphi)\cong\ZZ/2} -- is also
non-zero,
which implies that $\varphi$ cannot be realized by a map of spaces
\w{f\co \bX\to\bS{n-1}} (or even of suitable Postnikov sections).

\begin{remark}\label{rnonreal}
There is a more elementary way to see that $\varphi$ is not realizable:
if it were, from \eqref{eqtwentynine} and \eqref{eqthirty} we would have
\begin{multline}
\{6\nu,18\nu\}=\{6\nu,6\nu+\eta^{3}\}=
\varphi\{\beta,\beta+\alpha\circ\eta^{2}\}=f_{\ast}(\lra{\eta,2,\alpha})\\
=
\lra{\eta,2,\varphi(\alpha)}=
\lra{\eta,2,\eta}=\{\nu,12\nu\},
\end{multline}
\noindent a contradiction. Nevertheless, we hope the cohomological
approach helps to illustrate how the general theory works.
\end{remark}
\pagebreak
\bibliographystyle{gtart}
\bibliography{link}

\begin{thebibliography}{}
\providecommand\bibmarginpar{\leavevmode\marginpar}
\def\urlstyle#1{{\tt #1}}

\bibitem{ARosiL}
\textbf{J Ad{\'a}mek}, \textbf{J Rosick{\'y}}, \emph{Locally presentable and
  accessible categories}, London Mathematical Society Lecture Note Series 189,
  Cambridge University Press, Cambridge (1994) \xox{MR}{1294136}

\bibitem{BauCF}
\textbf{H-J Baues}, \emph{Combinatorial foundation of homology and homotopy},
  Springer Monographs in Mathematics, Springer, Berlin (1999) \xox{MR}{1707308}

\bibitem{BKS}
\textbf{D Benson}, \textbf{H Krause}, \textbf{S Schwede},
  \href{http://dx.doi.org/10.1090/S0002-9947-03-03373-7} {\emph{Realizability
  of modules over {T}ate cohomology}}, Trans. Amer. Math. Soc. 356 (2004)
  3621--3668 \xox{MR}{2055748}

\bibitem{BlaCH}
\textbf{D Blanc}, \emph{Comparing homotopy categories}, K--Theory (to appear)

\bibitem{BlaH}
\textbf{D\,A Blanc},
  \href{http://links.jstor.org/sici?sici=0002-9947(199003)318:1%3C335:AHSSFH%3%
E2.0.CO%3B2--P} {\emph{A {H}urewicz spectral sequence for homology}}, Trans.
  Amer. Math. Soc. 318 (1990) 335--354 \xox{MR}{956029}

\bibitem{BlaHH}
\textbf{D Blanc}, \emph{Higher homotopy operations and the realizability of
  homotopy groups}, Proc. London Math. Soc. $(3)$ 70 (1995) 214--240
  \xox{MR}{1300845}

\bibitem{BlaM}
\textbf{D Blanc}, \emph{Mapping spaces and $M$--CW complexes}, Forum Math. 9
  (1997) 367--382 \xox{MR}{1441926}

\bibitem{BlaAI}
\textbf{D Blanc}, \href{http://dx.doi.org/10.1017/S030500419900393X}
  {\emph{Algebraic invariants for homotopy types}}, Math. Proc. Cambridge
  Philos. Soc. 127 (1999) 497--523 \xox{MR}{1713125}

\bibitem{BlaCW}
\textbf{D Blanc}, \href{http://dx.doi.org/10.1016/S0166-8641(98)00090--X}
  {\emph{CW simplicial resolutions of spaces with an application to loop
  spaces}}, Topology Appl. 100 (2000) 151--175 \xox{MR}{1733041}

\bibitem{BDGoerR}
\textbf{D Blanc}, \textbf{W\,G Dwyer}, \textbf{P\,G Goerss},
  \href{http://dx.doi.org/10.1016/S0040-9383(03)00074-0} {\emph{The realization
  space of a $\Pi$--algebra: a moduli problem in algebraic topology}}, Topology
  43 (2004) 857--892 \xox{MR}{2061210}

\bibitem{BPescF}
\textbf{D Blanc}, \textbf{G Peschke}, \emph{The fiber of functors between
  categories of algebras}, J. Pure. Appl. Alg. (to appear)

\bibitem{BouC}
\textbf{A\,K Bousfield}, \href{http://dx.doi.org/10.2140/gt.2003.7.1001}
  {\emph{Cosimplicial resolutions and homotopy spectral sequences in model
  categories}}, Geom. Topol. 7 (2003) 1001--1053 \xox{MR}{2026537}

\bibitem{BFrH}
\textbf{A\,K Bousfield}, \textbf{E\,M Friedlander}, \emph{Homotopy theory of
  $\Gamma $--spaces, spectra, and bisimplicial sets}, from: ``Geometric
  applications of homotopy theory (Proc. Conf., Evanston, Ill., 1977), II'',
  Lecture Notes in Math. 658, Springer, Berlin (1978)  80--130 \xox{MR}{513569}

\bibitem{CDIComp}
\textbf{W Chach{\'o}lski}, \textbf{W\,G Dwyer}, \textbf{M Intermont},
  \href{http://dx.doi.org/10.1112/S0024610701002691} {\emph{The $A$--complexity
  of a space}}, J. London Math. Soc. $(2)$ 65 (2002) 204--222 \xox{MR}{1875145}

\bibitem{DroC}
\textbf{E Dror~Farjoun}, \emph{Cellular inequalities}, from: ``The Cech
  centennial (Boston, 1993)'', Contemp. Math. 181, Amer. Math. Soc.,
  Providence, RI (1995)  159--181 \xox{MR}{1320991}

\bibitem{DKanC}
\textbf{W\,G Dwyer}, \textbf{D\,M Kan},
  \href{http://dx.doi.org/10.1016/0040-9383(84)90035-1} {\emph{A classification
  theorem for diagrams of simplicial sets}}, Topology 23 (1984) 139--155
  \xox{MR}{744846}

\bibitem{DKanO}
\textbf{W\,G Dwyer}, \textbf{D\,M Kan},
  \href{http://dx.doi.org/10.1016/0040-9383(84)90035-1} {\emph{An obstruction
  theory for diagrams of simplicial sets}}, Nederl. Akad. Wetensch. Indag.
  Math. 46 (1984) 139--146 \xox{MR}{749527}

\bibitem{DKanR}
\textbf{W\,G Dwyer}, \textbf{D\,M Kan},
  \href{http://links.jstor.org/sici?sici=0002-9939(198407)91:3%3C456:RDITHC%3E%
2.0.CO%3B2--R} {\emph{Realizing diagrams in the homotopy category by means of
  diagrams of simplicial sets}}, Proc. Amer. Math. Soc. 91 (1984) 456--460
  \xox{MR}{744648}

\bibitem{DKanCR}
\textbf{W\,G Dwyer}, \textbf{D\,M Kan},
  \href{http://links.jstor.org/sici?sici=0002-9939(199202)114:2%3C575:CMAROD%3%
E2.0.CO%3B2--H} {\emph{Centric maps and realization of diagrams in the homotopy
  category}}, Proc. Amer. Math. Soc. 114 (1992) 575--584 \xox{MR}{1070515}

\bibitem{DKSmiR}
\textbf{W\,G Dwyer}, \textbf{D\,M Kan}, \textbf{J\,H Smith},
  \href{http://dx.doi.org/10.1016/0022-4049(89)90023-6} {\emph{Homotopy
  commutative diagrams and their realizations}}, J. Pure Appl. Algebra 57
  (1989) 5--24 \xox{MR}{984042}

\bibitem{DKStoE}
\textbf{W\,G Dwyer}, \textbf{D\,M Kan}, \textbf{C\,R Stover},
  \href{http://dx.doi.org/10.1016/0022-4049(93)90126--E} {\emph{An $E^2$ model
  category structure for pointed simplicial spaces}}, J. Pure Appl. Algebra 90
  (1993) 137--152 \xox{MR}{1250765}

\bibitem{DKStoB}
\textbf{W\,G Dwyer}, \textbf{D\,M Kan}, \textbf{C\,R Stover},
  \href{http://dx.doi.org/10.1016/0022-4049(94)00102--O} {\emph{The bigraded
  homotopy groups $\pi_{i,j}X$ of a pointed simplicial space $X$}}, J. Pure
  Appl. Algebra 103 (1995) 167--188 \xox{MR}{1358761}

\bibitem{ESteH}
\textbf{G Ellis}, \textbf{R Steiner},
  \href{http://dx.doi.org/10.1016/0022-4049(87)90089-2}
  {\emph{Higher-dimensional crossed modules and the homotopy groups of
  $(n{+}1)$--ads}}, J. Pure Appl. Algebra 46 (1987) 117--136 \xox{MR}{897011}

\bibitem{GHopkR}
\textbf{P\,G Goerss}, \textbf{M\,J Hopkins}, \emph{Resolutions in model
  categories}, preprint (1999)

\bibitem{GHopkM}
\textbf{P\,G Goerss}, \textbf{M\,J Hopkins}, \emph{Moduli spaces of commutative
  ring spectra}, from: ``Structured ring spectra'', London Math. Soc. Lecture
  Notes 315, Cambridge Univ. Press, Cambridge (2004)  151--200
  \xox{MR}{2125040}

\bibitem{GHopkM2}
\textbf{P\,G Goerss}, \textbf{M\,J Hopkins}, \emph{Moduli problems for
  structured ring spectra}, preprint (2005)
\ Available at \setbox0\hbox{\makeatletter\@url
{http://hopf.math.purdue.edu/cgi-bin/generate?/\allowbreak Goerss-Hopkins/obstruct}}
\href{http://hopf.math.purdue.edu/cgi-bin/generate?/Goerss-Hopkins/obstruct}
{\unhbox0}

\bibitem{GJarS}
\textbf{P\,G Goerss}, \textbf{J\,F Jardine}, \emph{Simplicial homotopy theory},
  Progress in Mathematics 174, Birkh\"auser Verlag, Basel (1999)
  \xox{MR}{1711612}

\bibitem{HirM}
\textbf{P\,S Hirschhorn}, \emph{Model categories and their localizations},
  Mathematical Surveys and Monographs 99, American Mathematical Society,
  Providence, RI (2003) \xox{MR}{1944041}

\bibitem{IlC}
\textbf{L Illusie}, \emph{Complexe cotangent et d\'eformations I}, Springer,
  Berlin (1971) \xox{MR}{0491680}

\bibitem{JardB}
\textbf{J\,F Jardine}, \emph{Bousfield's $E_2$ model theory for simplicial
  objects}, from: ``Homotopy theory: relations with algebraic geometry, group
  cohomology, and algebraic $K$--theory'', Contemp. Math. 346, Amer. Math.
  Soc., Providence, RI (2004)  305--319 \xox{MR}{2066504}

\bibitem{LodaS}
\textbf{J-L Loday}, \href{http://dx.doi.org/10.1016/0022-4049(82)90014-7}
  {\emph{Spaces with finitely many nontrivial homotopy groups}}, J. Pure Appl.
  Algebra 24 (1982) 179--202 \xox{MR}{651845}

\bibitem{MMSShipM}
\textbf{M\,A Mandell}, \textbf{J\,P May}, \textbf{S Schwede}, \textbf{B
  Shipley}, \href{http://dx.doi.org/10.1112/S0024611501012692} {\emph{Model
  categories of diagram spectra}}, Proc. London Math. Soc. $(3)$ 82 (2001)
  441--512 \xox{MR}{1806878}

\bibitem{MayS}
\textbf{J\,P May}, \emph{Simplicial objects in algebraic topology}, Van
  Nostrand Mathematical Studies 11, D. Van Nostrand Co., Princeton,
  N.J.-Toronto, Ont.-London (1967) \xox{MR}{0222892}

\bibitem{QuiS}
\textbf{D\,G Quillen}, \emph{Spectral sequences of a double semi-simplicial
  group}, Topology 5 (1966) 155--157 \xox{MR}{0195097}

\bibitem{QuiH}
\textbf{D\,G Quillen}, \emph{Homotopical algebra}, Lecture Notes in Mathematics
  43, Springer, Berlin (1967) \xox{MR}{0223432}

\bibitem{SpalSM}
\textbf{J Spali{\'n}ski}, \emph{Stratified model categories}, Fund. Math. 178
  (2003) 217--236 \xox{MR}{2030483}

\bibitem{StoV}
\textbf{C\,R Stover}, \href{http://dx.doi.org/10.1016/0040-9383(90)90022--C}
  {\emph{A van {K}ampen spectral sequence for higher homotopy groups}},
  Topology 29 (1990) 9--26 \xox{MR}{1046622}

\bibitem{TodaC}
\textbf{H Toda}, \emph{Composition methods in homotopy groups of spheres},
  Annals of Mathematics Studies 49, Princeton University Press, Princeton, N.J.
  (1962) \xox{MR}{0143217}

\bibitem{GWhE}
\textbf{G\,W Whitehead}, \emph{Elements of homotopy theory}, Graduate Texts in
  Mathematics 61, Springer, New York (1978) \xox{MR}{516508}

\end{thebibliography}

\end{document}